\documentclass[10pt,a4paper]{article}

\usepackage[a4paper,width=125mm,height=195mm,footskip=10mm]{geometry}
\usepackage[utf8]{inputenc}
\usepackage[T1]{fontenc}
\usepackage{lmodern}
\usepackage[english]{babel}
\usepackage{csquotes}
\usepackage{calc}
\usepackage{authblk}
\usepackage{fancyhdr}
\usepackage{titlesec}
\titlespacing*{\section}{0pt}{3.2ex plus 1ex minus .2ex}{1.8ex plus .2ex}
\titlespacing*{\subsection}{0pt}{3.0ex plus 1ex minus .2ex}{1.2ex plus .2ex}
\titlespacing*{\subsubsection}{0pt}{3.0ex plus 1ex minus .2ex}{1.2ex plus .2ex}
\titlespacing*{\paragraph}{0pt}{1.2ex plus 1ex minus .2ex}{1em}
\titlespacing*{\subparagraph}{\parindent}{1.2ex plus 1ex minus .2ex}{1em}

\usepackage{amsmath}
\usepackage{amsfonts}
\usepackage{amssymb}
\usepackage{mathrsfs}
\usepackage{amsthm}
\usepackage{thmtools}
\usepackage{mathtools}
\usepackage{bbm}

\allowdisplaybreaks

\mathtoolsset{showonlyrefs}

\let\originalleft\left
\let\originalright\right
\renewcommand{\left}{\mathopen{}\mathclose\bgroup\originalleft}
\renewcommand{\right}{\aftergroup\egroup\originalright}

\theoremstyle{plain}
\newtheorem{theorem}{Theorem}
\newtheorem*{theorem*}{Theorem}

\newtheorem{lemma}{Lemma}

\theoremstyle{definition}
\newtheorem{definition}{Definition}

\newtheorem*{remarkxx}{Remark}

\newtheorem*{examplexx}{Example}

\newenvironment{remark*}
{\pushQED{\qed}\remarkxx}
{\popQED\endremarkxx\vspace{2ex}}

\newenvironment{example*}
{\pushQED{\qed}\examplexx}
{\popQED\endexamplexx\vspace{2ex}}

\usepackage{float}
\usepackage{caption}
\usepackage{subcaption}
\usepackage{color}
\usepackage[table]{xcolor}
\usepackage{tabularx}
\usepackage{booktabs}
\usepackage{tikz}
\usepackage{pgfplots}

\usetikzlibrary{calc,decorations.pathmorphing,decorations.markings,shapes,arrows.meta,fadings}
\pgfplotsset{compat=1.12}

\usepackage{enumitem}

\setlist{align=right,labelindent=0.3em,labelwidth=1.6em,labelsep*=0.6em,leftmargin=!,topsep=0.5ex,partopsep=0.5ex,parsep=0ex,itemsep=0.25ex}
\setlist[itemize,1]{label=\raisebox{0.0em}{\scalebox{1.0}{$\bullet$}}}
\setlist[itemize,2]{label=\raisebox{0.1em}{\scalebox{0.5}{$\blacksquare$}}}
\setlist[itemize,3]{label=\raisebox{0.1em}{\scalebox{0.7}{$\blacktriangleright$}}}
\setlist[itemize,4]{label=\raisebox{0.1em}{\scalebox{0.6}{$\blacklozenge$}}}
\setlist[enumerate]{label=(\alph*)}

\usepackage[backend=biber,style=numeric,sorting=nyt,sortcites=true,doi=false,isbn=false,url=false,eprint=true,maxbibnames=9,giveninits=true]{biblatex}
\usepackage[colorlinks=false,bookmarksnumbered=false,pdfpagemode=UseNone,pdfstartpage=1,pdfstartview={XYZ null null 1.00},pdfpagelayout=OneColumn,pdflang=English]{hyperref}

\DeclareNameAlias{sortname}{last-first}
\DeclareNameAlias{author}{first-last}

\renewbibmacro{in:}{}

\AtEveryBibitem{\ifentrytype{misc}{}{%
	\clearfield{archivePrefix}%
	\clearfield{arxivId}%
	\clearfield{eprint}%
	\clearfield{primaryClass}%
}}

\usepackage{todonotes}

\renewcommand{\P}{\ensuremath{\mathbb{P}}}
\newcommand{\E}{\ensuremath{\mathbb{E}}}

\newcommand{\e}{\ensuremath{\mathrm{e}}}
\renewcommand{\epsilon}{\varepsilon}

\newcommand{\1}[1]{\ensuremath{\mathbbm{1}{\raisebox{-0.14ex}{\hspace{-0.06em}{}$\scriptstyle #1$}}}} % indicator function
\newcommand{\symdiff}{\ensuremath{\mathbin{\scalebox{0.9}{$\bigtriangleup$}}}} %symmetric difference of sets
\newcommand{\dif}{\mathop{}\!\mathrm{d}} % differential operator (the 'd' in integrals)
\newcommand{\conv}[1][\P]{\mathrel{\xrightarrow{\scalebox{0.65}{\raisebox{-0.5pt}[0pt][0pt]{$\scriptstyle{}#1$}}}}}

\DeclareMathOperator{\lambertW}{\scalebox{0.85}{$\mathcal{W}$}}

\newcommand{\bigO}{\ensuremath{\mathord{O}}}

\newcommand{\smallO}{\ensuremath{\mathord{o}}}

\newcommand{\bigTheta}{\ensuremath{\mathord{\Theta}}}

\newcommand{\bigOmega}{\ensuremath{\mathord{\Omega}}}

\newcommand{\grg}{\ensuremath{\mathbb{G}(n, d, p)}}    % Geometric random graph (null hypothesis)
\newcommand{\brg}{\ensuremath{\mathbb{G}(n, d, p; k)}} % Botnet random graph    (alt. hypothesis)

\newcommand{\torus}{\textup{T}}
\newcommand{\distG}{D_{\textup{\tiny{}G}\mspace{-1mu}}}
\newcommand{\distT}{D_{\textup{\tiny{}T}\mspace{-1mu}}}
\newcommand{\diam}{\textup{diam}}

\newcommand{\test}{\psi}

\newcommand{\est}{\ensuremath{\widehat{B}}}
\newcommand{\risk}{\ensuremath{R}}
\newcommand{\detrisk}{\ensuremath{\risk}}
\newcommand{\estrisk}{\ensuremath{\risk_{\textup{est}}}}

\newcommand{\ball}{\ensuremath{B}}
\newcommand{\kissingnumber}[1]{\ensuremath{\kappa_{#1}}}
\newcommand{\ballsubset}{A}

\newcommand{\partdeg}{\text{deg}_{\scriptscriptstyle{}V \setminus B}}

\makeatletter %% <- make @ usable in command names

\newcommand*\sNeg[2][0mu]{\Neginternal{#1}{\snegslash}{#2}}

\newcommand*\Neginternal[3]{\mathpalette\Neg@{{#1}{#2}{#3}}}
\newcommand*\Neg@[2]{\Neg@@{#1}#2}
\newcommand*\Neg@@[4]{%
  \mathrel{\ooalign{%
    $\m@th#1#4$\cr
    \hidewidth$\m@th#3{#1}\mkern\muexpr#2*2$\hidewidth\cr
  }}%
}
\newcommand*\negslash[1]{\m@th#1\not\mathrel{\phantom{=}}}
\newcommand*\snegslash[1]{\rotatebox[origin=c]{60}{$\m@th#1-$}}
\newcommand*\ssnegslash[1]{\rotatebox[origin=c]{60}{$\m@th#1{\dabar@}\mkern-7mu{\dabar@}$}}
\newcommand*\sssnegslash[1]{\rotatebox[origin=c]{60}{$\m@th#1\dabar@$}}

\newcounter{sarrow}
\newcommand{\@leftrightsquigarrow}[1][\rule{6.5pt}{0pt}]{\stepcounter{sarrow}
\mathrel{\begin{tikzpicture}[baseline= {( $ (current bounding box.south) $ )}]
\node[inner sep=.5ex] (\thesarrow) {$\scriptstyle\smash{#1}$};
\path[draw,To-To,decorate,decoration={zigzag,amplitude=0.7pt,segment length=2.2pt,pre=lineto,pre length=2.4pt,post=lineto,post length=2pt}] (\thesarrow.south west) -- (\thesarrow.south east);
\end{tikzpicture}}}
\renewcommand{\leftrightsquigarrow}{\mathrel{\mathchoice{\@leftrightsquigarrow}{\@leftrightsquigarrow}{\scalebox{0.75}{$\mspace{2mu}\@leftrightsquigarrow\mspace{2mu}$}}{\scalebox{0.6}{$\mspace{1mu}\@leftrightsquigarrow\mspace{1mu}$}}}}
\makeatother

\renewcommand{\nleftrightarrow}{\mathrel{\sNeg\leftrightarrow}}

\newcommand{\eventib}{\mathord{\mathsf{A}}}
\newcommand{\eventcg}{\mathord{\mathsf{C}}}
\newcommand{\eventmd}{\mathord{\mathsf{D}}}

\title{Detecting a botnet in a network}
\author[{}\hspace{0.5pt}\protect\hyperlink{hyp:email1}{1},\protect\hyperlink{hyp:affil2}{b}]{\protect\hypertarget{hyp:author1}{Gianmarco Bet}}
\author[{}\hspace{0.5pt}\protect\hyperlink{hyp:email2}{2},\protect\hyperlink{hyp:affil1}{a},\protect\hyperlink{hyp:corresponding}{$\dagger$}]{\protect\hypertarget{hyp:author2}{Kay Bogerd}}
\author[{}\hspace{0.5pt}\protect\hyperlink{hyp:email3}{3},\protect\hyperlink{hyp:affil1}{a}]{\protect\hypertarget{hyp:author3}{Rui M. Castro}}
\author[{}\hspace{0.5pt}\protect\hyperlink{hyp:email4}{4},\protect\hyperlink{hyp:affil1}{a}]{\protect\hypertarget{hyp:author4}{Remco van der Hofstad}}
\affil[ ]{\small\parbox{365pt}{\parbox{5pt}{\textsuperscript{\protect\hypertarget{hyp:affil1}{a}}}Eindhoven University of Technology,%
\enspace\parbox{5pt}{\textsuperscript{\protect\hypertarget{hyp:affil2}{b}}}Università degli Studi di Firenze,}}
\affil[ ]{\small\parbox{365pt}{
\parbox{5pt}{\textsuperscript{\protect\hypertarget{hyp:email1}{1}}}\texttt{\footnotesize\href{mailto:gianmarco.bet@unifi.it}{gianmarco.bet@unifi.it}},
\parbox{5pt}{\textsuperscript{\protect\hypertarget{hyp:email2}{2}}}\texttt{\footnotesize\href{mailto:k.m.bogerd@tue.nl}{k.m.bogerd@tue.nl}},
\parbox{5pt}{\textsuperscript{\protect\hypertarget{hyp:email3}{3}}}\texttt{\footnotesize\href{mailto:rmcastro@tue.nl}{rmcastro@tue.nl}},
\parbox{5pt}{\textsuperscript{\protect\hypertarget{hyp:email4}{4}}}\texttt{\footnotesize\href{mailto:r.w.v.d.hofstad@tue.nl}{r.w.v.d.hofstad@tue.nl}},}}
\affil[ ]{\small\parbox{365pt}{\parbox{5pt}{\textsuperscript{\protect\hypertarget{hyp:corresponding}{$\dagger$}}}Corresponding author}}
\date{\today}\footnotesize

\begin{document}
\vspace*{-30pt}
\noindent\makebox[\textwidth][c]{\begin{minipage}[c]{1.25\textwidth}\maketitle\end{minipage}}

\begin{abstract}
We formalize the problem of detecting the presence of a botnet in a network as a hypothesis testing problem where we observe a single instance of a graph. The null hypothesis, corresponding to the absence of a botnet, is modeled as a random geometric graph where every vertex is assigned a location on a $d$-dimensional torus and two vertices are connected when their distance is smaller than a certain threshold. The alternative hypothesis is similar, except that there is a small number of vertices, called the botnet, that ignore this geometric structure and simply connect randomly to every other vertex with a prescribed probability.

We present two tests that are able to detect the presence of such a botnet. The first test is based on the idea that botnet vertices tend to form large isolated stars that are not present under the null hypothesis. The second test uses the average graph distance, which becomes significantly shorter under the alternative hypothesis. We show that both these tests are asymptotically optimal. However, numerical simulations show that the isolated star test performs significantly better than the average distance test on networks of moderate size. Finally, we construct a robust scheme based on the isolated star test that is also able to identify the vertices in the botnet.
\end{abstract}

%%%%%%%%%%%%%%%%%%%%%%%%%%%%%%%%%%%%%%%%%%%%%%%%%%%%%%%%%%%%%%%%%%%%%%%%%%%%%%%%
%%%%%%%%%%%%%%%%%%%%%%%%%%%%%%%%%%%%%%%%%%%%%%%%%%%%%%%%%%%%%%%%%%%%%%%%%%%%%%%%
\section{Introduction}
\label{sec:introduction}
Complex networks are often described in terms of a large number of vertices that are connected using the same underlying probabilistic mechanism. In practice, however, these networks might contain a small number of vertices that follow different connection criteria. Examples are fake user profiles in a social network (like Facebook or LinkedIn) or servers infected by a computer virus on the internet. We refer to such a set of anomalous vertices as a \emph{botnet}. Typically a botnet represents a potentially malicious anomaly in the network, and thus it is of great practical interest to detect its presence and, when detected, to identify the corresponding vertices. Accordingly, numerous empirical studies have analyzed botnet detection problems and techniques, see \cite{Feily2009,Zeidanloo2010,Garcia2014a,Garcia2014,Mesnards2018} and the references therein. In this work we look at the problem from a statistical point of view, and characterize the difficulty of detecting a botnet based only on structural information from the observed network.

More precisely, we formalize this problem as a hypothesis testing problem where we observe a single instance of a random graph. Under the null hypothesis, this graph is a sample from a random geometric graph \cite{Gilbert1961,Penrose2003} on $n$ vertices where every vertex is assigned a location on a $d$-dimensional torus and two vertices are connected when their Euclidean distance on the torus is less than a given radius. Under the alternative hypothesis there is a small number $k$ of vertices, called the botnet, that ignore the geometric structure and instead connect to every other vertex with a prescribed probability. In other words, $n-k$ vertices still connect based on the underlying geometry, while each of the $k$ botnet vertices forms connections uniformly at random with every other vertex (botnet or not). In practice, botnets are built to imitate regular nodes in the network, and so we assume that the expected degree of every vertex is the same under the null and alternative hypothesis. This assumption rules out trivial scenarios where the botnet can be detected simply by looking at the edge density or degree structure.

\paragraph{Our contribution.} We propose two different tests to detect whether an observed graph contains a botnet. The first test is a local test, based on the number of isolated stars that can be observed in the given graph. For convenience we refer to this test as the \emph{isolated star test}. For a given vertex, its isolated star is the largest subset of its neighbors such that none of them are connected to each other by an edge. Hence, an isolated star is the largest independent set on the subgraph induced by the neighbors of a vertex. Under the null hypothesis, none of the vertices can become a large isolated star because the underlying geometry ensures that most neighbors are directly connected. However, because the botnet vertices are connected uniformly at random throughout the graph they are likely to become large isolated stars.

Our second test is based on graph distances in the observed graph and thus it has a more global nature. We refer to this test as the \emph{average distance test}. Under the null hypothesis, vertices that are separated by a large Euclidean distance will also be separated by a large graph distance. However, under the alternative hypothesis, the botnet vertices typically create shortcuts, making many paths much shorter. Under appropriate assumptions, the effect of the shortcuts is large enough to significantly decrease the average graph distance. This phenomenon was first investigated by Watts and Strogatz \cite{Watts1998}.

Both of our methods can be used to test for the presence of a botnet. Our results show that a botnet can be detected, with high probability, when the expected number of edges connected to all botnet vertices is diverging (i.e., when the expected vertex degree diverges or when the botnet size is unbounded). Remarkably, this means that a single botnet vertex can be detected provided that the graph is not of bounded average degree. We also show that this result is optimal, meaning that it is impossible for any test to detect the presence of a botnet when the expected number of botnet edges is bounded. We complement our theoretical results for the $n \to \infty$ asymptotic regime with numerical simulations that illustrate the performance of our tests on graphs of finite size. These results empirically show that the isolated star test performs much better than the average distance test, with the difference being more pronounced when the dimension of the underlying geometry is large.

\paragraph{Related work.} Recently there has been an increasing interest in the development of statistical techniques and algorithms that exploit the structure of large complex-network data to analyze networks more efficiently. In particular, several recent papers have studied hypothesis testing for random graph models. In \cite{Arias-Castro2015,Arias-Castro2014}, the authors consider the problem of detecting a denser subset of vertices in an Erd\H{o}s-R\'{e}nyi random graph, or in an inhomogeneous random graph \cite{Bogerd2019}.

The setting of \cite{Bubeck2016} is perhaps the closest to our setting. The authors consider the problem of deciding whether a given graph is generated by some underlying spatial mechanism. More specifically, in their model, the null hypothesis is an Erd\H{o}s-R\'{e}nyi random graph, and this is compared to a high-dimensional random geometric graph under the alternative. As the dimension tends to infinity, the two random graphs become indistinguishable, and they identify how large the dimension can be so that these models can still be distinguished.

The authors of \cite{Gao2017b} propose a test based on observed frequencies of small subgraphs to distinguish between an Erd\H{o}s-R\'{e}nyi random graph, seen as the null hypothesis, and a general class of alternative models that include stochastic block models and the configuration model. Similarly, \cite{Bresler2018} proposes a test to distinguish between mean-field models and structured Gibbs models. Finally, \cite{Heard2010,Mongiovi2013,Park2013} investigate detection problems in a dynamical setting, where the goal is to detect changes in the graph structure over time.

In this paper we specifically consider the problem of detecting a botnet in undirected graphs. For instance, servers infected by a computer virus on the internet or fake user profiles in social networks like Facebook or LinkedIn. A related and very interesting problem is that of detecting botnets in directed networks, such as Twitter. These are heavily involved in the spread of fake news \cite{Shah2011,Bhamidi2015,Mesnards2018,Crane2020}. In both settings we are trying to identify nodes in the network that are anomalous or disruptive. However, the way these anomalous nodes manifest themselves is rather different than in our model.

%%%%%%%%%%%%%%%%%%%%%%%%%%%%%%%%%%%%%%%%%%%%%%%%%%%%%%%%%%%%%%%%%%%%%%%%%%%%%%%%
%%%%%%%%%%%%%%%%%%%%%%%%%%%%%%%%%%%%%%%%%%%%%%%%%%%%%%%%%%%%%%%%%%%%%%%%%%%%%%%%
\section{Model formulation and results}
\label{sec:model_formulation_and_results}
In this section we formalize the problem of detecting a botnet in a network as a hypothesis testing problem for graphs. We are given a single observation of a random graph $G = (V, E)$, where $V = \{1, \ldots, n\}$ is the vertex set of size $|V| = n$ and $E \subseteq \{(i,j) \in V \times V : i < j\}$ is the random set of edges. We use $i \leftrightarrow j$ to indicate that $i,j \in V$ are connected. That is, we write $i \leftrightarrow j$ when $(i,j) \in E$ and $i \nleftrightarrow j$ otherwise. In particular, $G$ is a simple graph, so it does not contain any self-loops or multiple edges.

Under the null hypothesis, denoted by $H_0$, the observed graph $G$ is a realization of a $d$-dimensional random geometric graph $\grg$ on $n$ vertices and with average edge probability $p$. Formally, let $\torus^d \coloneqq [0, 1]^d$ be the $d$-dimensional unit torus, with distance function
\begin{equation}
\label{eq:torus_distance}
\distT(x, y)
  = \sqrt{{\textstyle\sum_{j = 1}^d} \min\bigl(|x_j - y_j|, 1 - |x_j - y_j|\bigr)^2} \,,
\qquad\text{for } x, y \in \torus^d \,.
\end{equation}
This is simply the Euclidean distance on the unit (hyper-)cube with the ability to ``wrap around'' the boundaries. We refer to $\torus^d$ as the embedding space. For each vertex $i \in V$, let $X_i$ be a $d$-dimensional vector-valued random variable uniformly distributed on $\torus^d$. We denote the components of this random vector by $X_i = (X_{i,1}, \ldots, X_{i,d})$ and note that these components are independent uniform random variables on the unit interval $[0, 1]$.

For a given edge probability $p$, two vertices $i, j \in V$ are connected when $\distT(X_i, X_j) \leq r$, where $r$ is chosen such that the average edge probability equals $p$, that is $\P(\distT(X_i, X_j) \leq r) = p$. In other words, $r$ is such that the probability of a random point $X_i$ landing in a ball of radius $r$ is equal to $p$, which gives the explicit relation $p = (\sqrt{\pi} \, r)^d / \Gamma(d/2 + 1)$, where $\Gamma(\cdot)$ denotes the gamma function. Throughout the rest of this paper we assume that $p \to 0$ as $n \to \infty$, so the average degree is sub-linear in the graph size $n$. For further details on this model and many of its properties we refer the reader to \cite{Penrose2003}.

The alternative hypothesis, denoted by $H_1$, is similar except for a small subset of vertices called the botnet. These vertices ignore the geometric structure and simply connect to every other vertex independently with probability $p$. Formally, the observed graph under the alternative hypothesis is a realization from $\brg$, which is a random geometric graph on $n-k$ vertices together with a subset of vertices $B \subseteq V$ of size $|B| = k$, called the botnet. That is, each pair of vertices $i, j \in V \setminus B$ is connected precisely when $\distT(X_i, X_j) \leq r$. The remaining vertices in the botnet $B$ are connected independently and with probability $p$ to every other vertex in $V$. Note that, by construction, the expected number of edges under the alternative hypothesis is exactly the same as under the null hypothesis.

\begin{figure}[H]
\centering
\begin{subfigure}[b]{0.49\textwidth}
\centering
\includegraphics[width=0.80\textwidth]{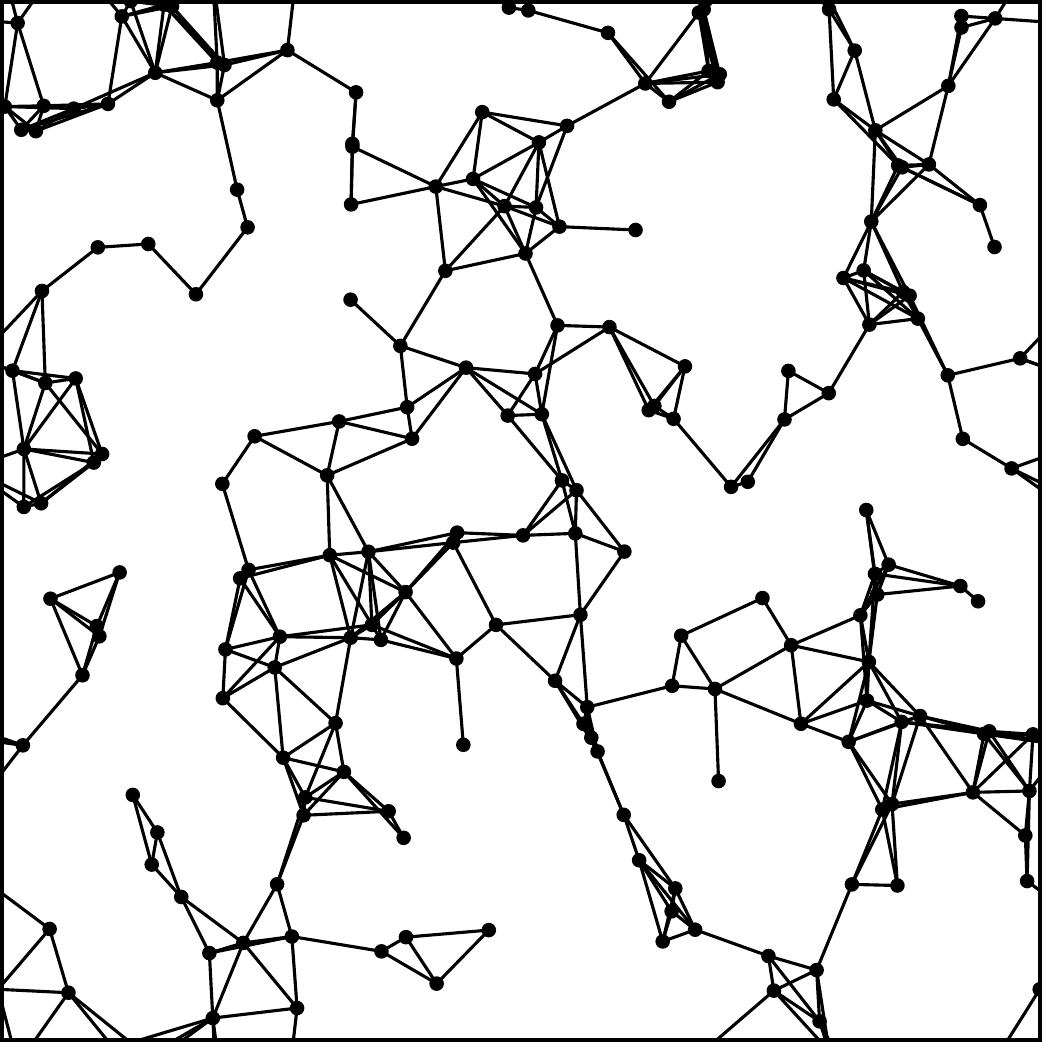}
\vspace{7pt}
\caption{Null model $\grg$.}
\label{fig:model_example_null}
\end{subfigure}%
\hspace{0.0199\textwidth}%
\begin{subfigure}[b]{0.49\textwidth}
\centering
\includegraphics[width=0.80\textwidth]{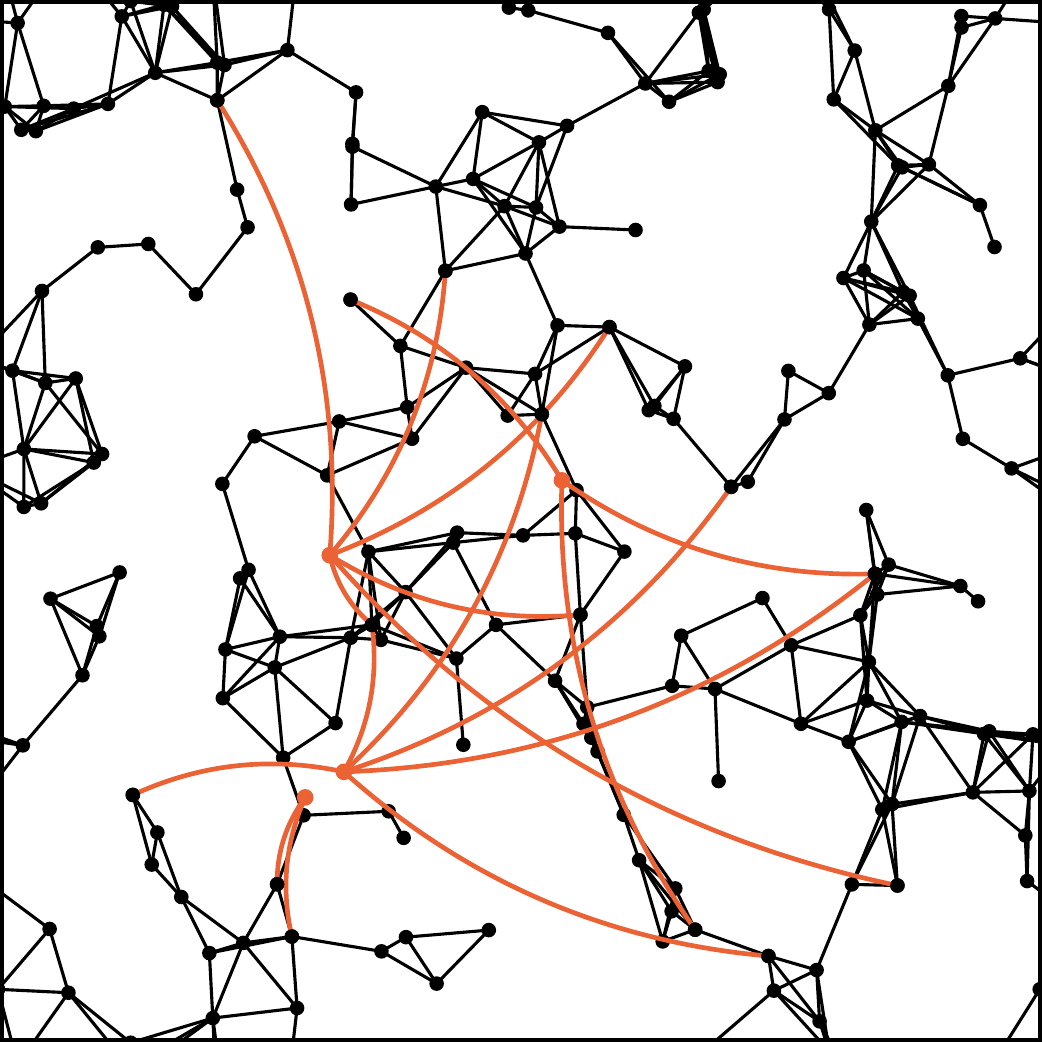}
\vspace{7pt}
\caption{Alternative model $\brg$.}
\label{fig:model_example_alt}
\end{subfigure}%
\vspace{4pt}
\caption{Example of the model under the null and alternative model in $2$ dimensions, were we identified opposite sides of the square so that edges can ``wrap-around'' the sides. Note that this representation uses the embedding of the vertices in the torus that is not available for the inference problem. The graph contains $n = 200$ vertices with $k = 4$ botnet vertices and average degree $n p = 5$. The botnet is highlighted in red.}
\label{fig:model_example}
\end{figure}

\begin{figure}[H]
\centering
\begin{subfigure}[b]{0.49\textwidth}
\centering
\vspace{-6pt}
\includegraphics[width=0.90\textwidth,angle=90]{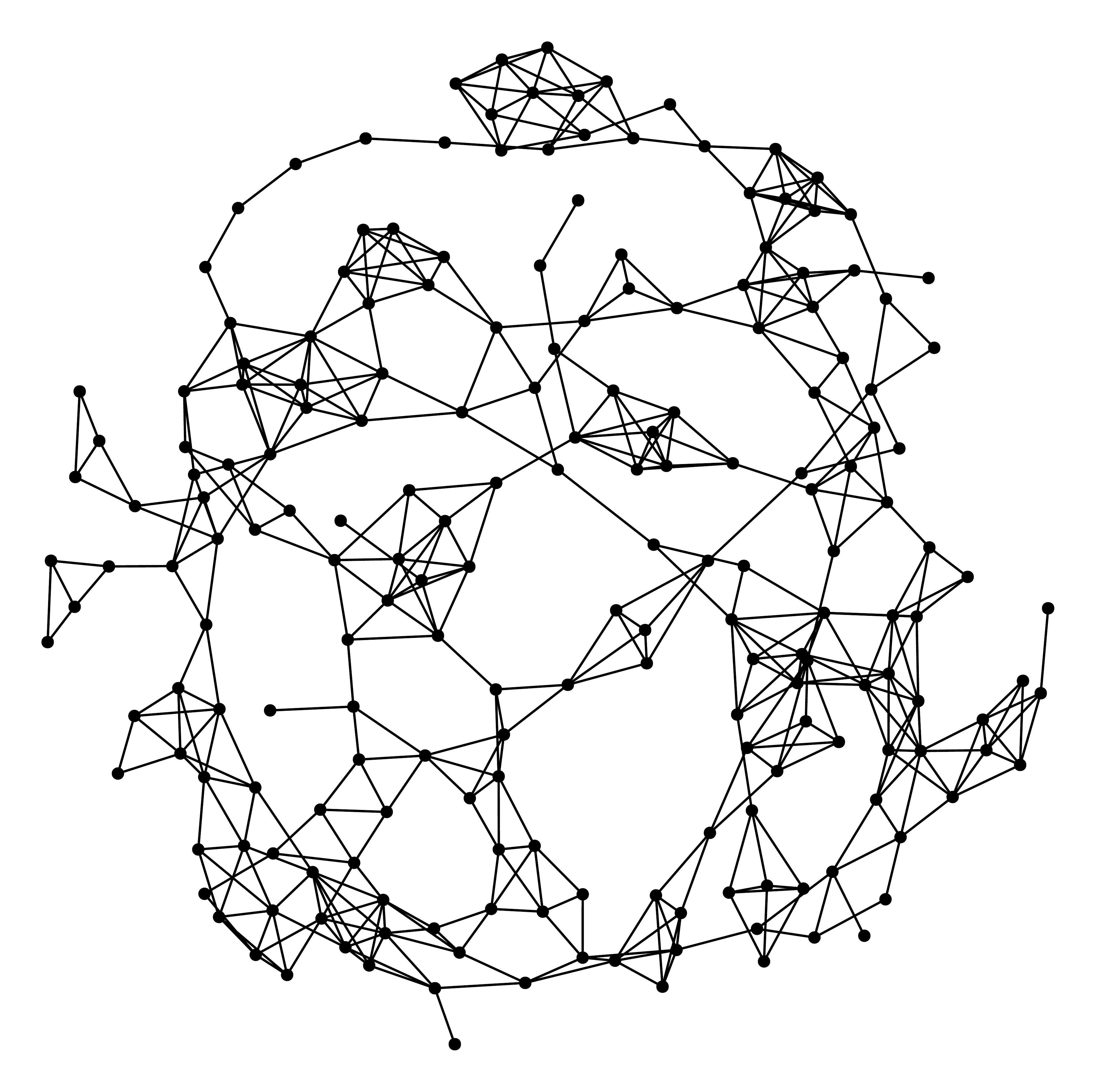}
\vspace{-6pt}
\caption{Null model $\grg$.}
\label{fig:model_example_null2}
\end{subfigure}%
\hspace{0.0199\textwidth}%
\begin{subfigure}[b]{0.49\textwidth}
\centering
\vspace{-6pt}
\includegraphics[width=0.90\textwidth,angle=90]{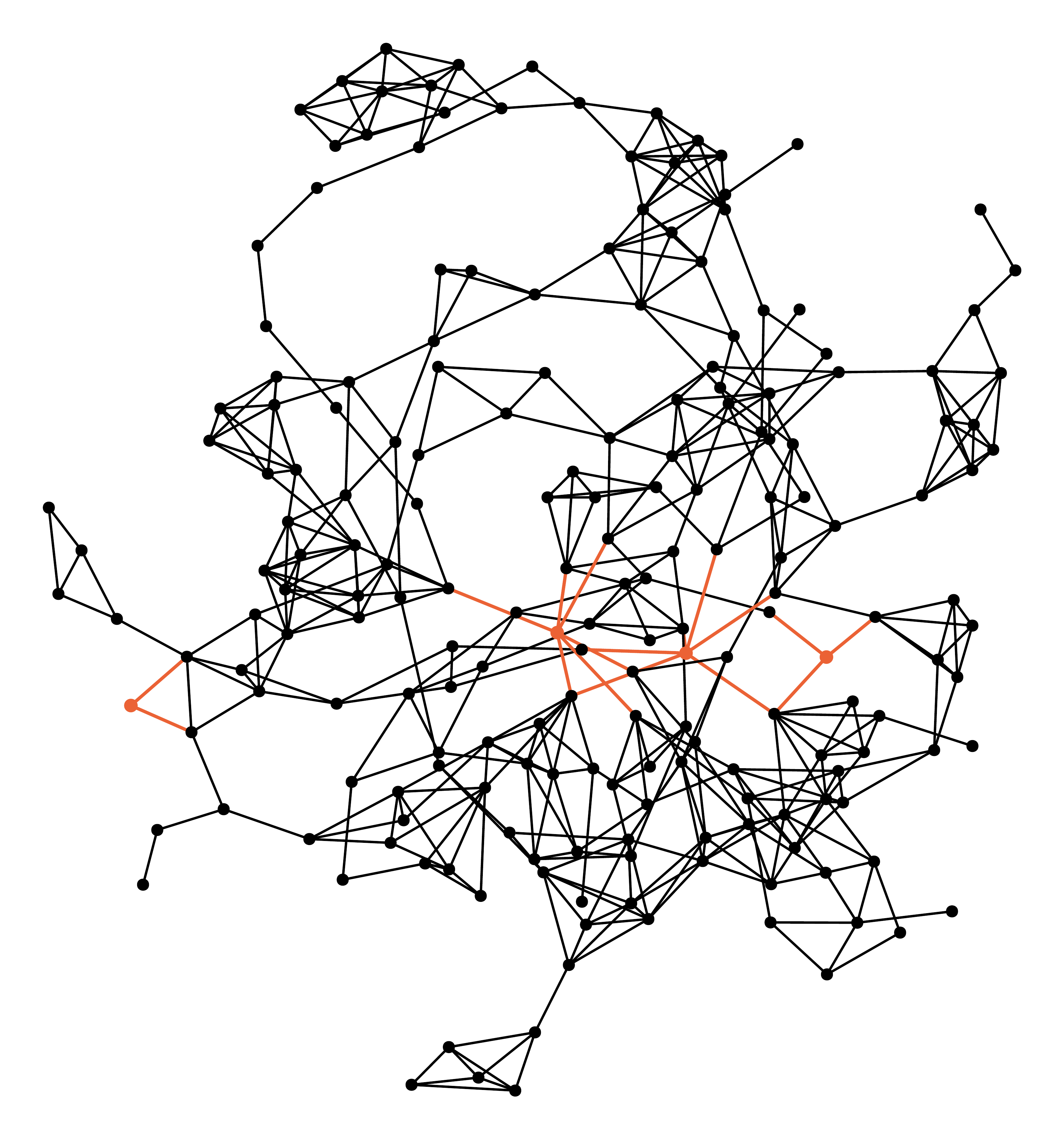}
\vspace{-6pt}
\caption{Alternative model $\brg$.}
\label{fig:model_example_alt2}
\end{subfigure}%
\vspace{2pt}
\caption{The same example graphs as in Figure~\mbox{\ref{fig:model_example}}, but drawn using a force field layout.}
\label{fig:model_example2}
\end{figure}

Another way to sample a graph $\brg$ from the alternative hypothesis is to first sample a graph $\grg$ from the null hypothesis. Then randomly select $k$ vertices and delete all edges incident to them, and finally reconnect these vertices to every other vertex independently and with probability $p$. An example of this is shown in Figures \ref{fig:model_example} and \ref{fig:model_example2}, where we compare the model under the null and alternative hypothesis in $2$ dimensions. However, remember that the vertex locations as shown in Figure~\ref{fig:model_example} are not available for the inference problem and we can only observe which vertices are connected. In Figure~\ref{fig:model_example2} a representation of the graph that does not rely on the Euclidean embedding is given, illustrating how the botnet edges faintly ``shorten'' the connections between different parts of the network.

\paragraph{General assumptions and notation.}
Throughout the rest of this paper all unspecified limits are assumed to be taken as the graph size $n$ tends to $\infty$. We also use standard asymptotic notation: $a_n = \bigO(b_n)$ when $a_n / b_n$ is bounded, $a_n = \bigOmega(b_n)$ when $b_n = \bigO(a_n)$, $a_n = \bigTheta(b_n)$ when $a_n = \bigO(b_n)$ and $a_n = \bigOmega(b_n)$, and $a_n = \smallO(b_n)$ when $a_n / b_n \to 0$. Furthermore, we write $a_n \asymp b_n$ to indicate that $a_n = (1 + \smallO(1)) b_n$, and $a_n \ll b_n$ when $a_n = \smallO(b_n)$, or $a_n \gg b_n$ when $b_n = \smallO(a_n)$. Finally, we say that a sequence of events holds with high probability if it holds with probability tending to $1$ as $n \to \infty$.

Given two vertices $i, j \in V$, we write $i \leftrightarrow j$ when these vertices are directly connected by an edge, and $i \leftrightsquigarrow j$ when there exists a path between them. Further, we assume that the dimension $d \geq 2$ remains fixed, but the edge probability $p$ and the botnet size $k$ are allowed to depend on $n$, although this dependence is left implicit in the notation. We also require that $p \to 0$ in such a way that $n p = \bigOmega(1)$ because otherwise the resulting graphs will be such that most vertices are isolated. Finally, we assume that the botnet size $k$ satisfies $1 \leq k \leq \smallO(n)$.

\vspace{10pt}
\vspace{-6pt}%

%%%%%%%%%%%%%%%%%%%%%%%%%%%%%%%%%%%%%%%%%%%%%%%%%%%%%%%%%%%%%%%%%%%%%%%%%%%%%%%%
%%%%%%%%%%%%%%%%%%%%%%%%%%%%%%%%%%%%%%%%%%%%%%%%%%%%%%%%%%%%%%%%%%%%%%%%%%%%%%%%
\subsection{Detecting a botnet}
\label{sec:detecting_a_botnet}
In this section we obtain a necessary condition for detecting the presence of a planted botnet in the asymptotic regime $n \to \infty$. Given an observed graph, we want to decide whether it was sampled from $H_0$ or from $H_1$. To this end, define a test $\test$ as a function mapping $G$ to $\{0, 1\}$, where $\test(G) = 1$ indicates the null hypothesis is rejected (i.e., the test indicates that the graph contains a botnet), and $\test(G) = 0$ otherwise. The \textit{worst-case risk} of such a test is defined as
\begin{equation}
\label{eq:worst_case_risk}
\detrisk(\test) \coloneqq \P_0(\test(G) \neq 0) + \max_{B \subseteq V,\, |B|=k} \P_B(\test(G) \neq 1) \,,
\end{equation}
where $\P_0(\cdot)$ denotes the distribution of the random geometric graph under the null hypothesis, and $\P_B(\cdot)$ denotes the distribution of a graph with the botnet $B \subseteq V$ under the alternative hypothesis.

Our goal is to determine when can we distinguish $H_0$ and $H_1$ as the graph size $n$ diverges. To this end we consider a sequence of tests $(\test_n)_{n=1}^{\infty}$ and we call such a sequence asymptotically powerful when it has vanishing risk, that is $\detrisk(\test_n) \to 0$ as $n \to \infty$. Hence, a sequence of tests is asymptotically powerful when it identifies the underlying model correctly in the limit $n \to \infty$.

Before we introduce our tests, we define the threshold (in terms of the model parameters) below which it becomes impossible for any test to be asymptotically powerful. We later show that above this threshold the isolated star test is asymptotically powerful. The average distance test is also asymptotically powerful in this regime, assuming some additional technical assumptions are satisfied. This threshold is given in terms of the parameters of the alternative model. Intuitively, it corresponds to the setting where the expected number of edges connected to all botnet vertices is bounded, which happens precisely when both the average degree $n p$ and the botnet size $k$ are bounded. In this case, there is a positive probability that all botnet vertices are isolated. When this happens it becomes impossible to reliably distinguish the null and alternative hypothesis. This is formalized in the following theorem, the proof of which is postponed to Section~\ref{sec:proof_no_test_powerful}:
\begin{theorem}
\label{thm:no_test_powerful}
When $n p k = \bigO(1)$ no test can be asymptotically powerful (i.e., all tests have risk that is strictly larger than zero).
\end{theorem}
In the rest of this section we present the two different tests that can detect the presence of a planted botnet in the regime $n p k \to \infty$. 

%%%%%%%%%%%%%%%%%%%%%%%%%%%%%%%%%%%%%%%%%%%%%%%%%%%%%%%%%%%%%%%%%%%%%%%%%%%%%%%%
%%%%%%%%%%%%%%%%%%%%%%%%%%%%%%%%%%%%%%%%%%%%%%%%%%%%%%%%%%%%%%%%%%%%%%%%%%%%%%%%
\subsubsection{Isolated star test}
\label{sec:isolated_star_test}
In this section we define a test that can detect whether an observed graph contains a planted botnet based on the presence of isolated stars. For a given vertex $i \in V$, let $N(i) =\{j \in V : (i,j) \in E\}$ denote the subset of its neighbors. The isolated star $S(i) \subseteq N(i)$, at vertex $i \in V$, is the largest independent set on the subgraph of $G$ induced by $N(i)$. In other words, every $j \in S(i)$ is directly connected by an edge to $i$, and no pair of vertices in $S(i)$ are directly connected (i.e., for every $j, k \in S(i)$ we have $(j,k) \notin E$).

\begin{figure}[thb]
\centering
\begin{tikzpicture}[scale=1.6]
\draw[white!80!black,line width=0.8pt,fill=white!95!black] (0.0,0.0) circle (1.0);

\draw[fill=black] (0.0,0.0) circle (0.03);
\draw (0.0,0.0) node[anchor=north west,inner xsep=-1.5pt,inner ysep=7pt,scale=0.68] {\Large$i$};

\draw[fill=black] (0.627979,-0.655302) circle (0.03);
\draw (0.627979,-0.655302) node[anchor=north east,inner xsep=3pt,inner ysep=-1.5pt,scale=0.6] {\Large$1$};
\draw (0,0) -- (0.627979,-0.655302);
\draw[path fading=north,fading transform={rotate=223.78}] (0.627979,-0.655302) -- (0.800952,-0.835802);
\draw[path fading=north,fading transform={rotate=183.78}] (0.627979,-0.655302) -- (0.644461,-0.904758);
\draw[path fading=north,fading transform={rotate=263.78}] (0.627979,-0.655302) -- (0.876507,-0.682387);

\draw[fill=black] (-0.327942,-0.834129) circle (0.03);
\draw (-0.327942,-0.834129) node[anchor=west,inner xsep=5pt,inner ysep=0pt,scale=0.6] {\Large$2$};
\draw (0,0) -- (-0.327942,-0.834129);
\draw[path fading=north,fading transform={rotate=106.037}] (-0.327942,-0.834129) -- (-0.568212,-0.933967);
\draw[path fading=north,fading transform={rotate=141.037}] (-0.327942,-0.834129) -- (-0.485145,-1.03519);
\draw[path fading=north,fading transform={rotate=176.037}] (-0.327942,-0.834129) -- (-0.345218,-1.08254);
\draw[path fading=north,fading transform={rotate=211.037}] (-0.327942,-0.834129) -- (-0.199042,-1.06333);

\draw[fill=black] (-0.964169,0.0586606) circle (0.03);
\draw (-0.964169,0.0586606) node[anchor=north west,inner xsep=1pt,inner ysep=4pt,scale=0.6] {\Large$3$};
\draw (0,0) -- (-0.964169,0.0586606);
\draw[path fading=north,fading transform={rotate=86.5184}] (-0.964169,0.0586606) -- (-1.21371,0.0738427);
\draw[path fading=north,fading transform={rotate=46.5184}] (-0.964169,0.0586606) -- (-1.14557,0.230691);
\draw[path fading=north,fading transform={rotate=126.518}] (-0.964169,0.0586606) -- (-1.16509,-0.0901095);

\draw[fill=black] (-0.341645,0.841691) circle (0.03);
\draw (-0.341645,0.841691) node[anchor=north east,inner xsep=4pt,inner ysep=1pt,scale=0.6] {\Large$4$};
\draw (0,0) -- (-0.341645,0.841691);
\draw[path fading=north,fading transform={rotate=329.592}] (-0.341645,0.841691) -- (-0.215108,1.07197);
\draw[path fading=north,fading transform={rotate=4.59239}] (-0.341645,0.841691) -- (-0.361661,1.08978);
\draw[path fading=north,fading transform={rotate=39.5924}] (-0.341645,0.841691) -- (-0.500975,1.04111);
\draw[path fading=north,fading transform={rotate=74.5924}] (-0.341645,0.841691) -- (-0.58266,0.939002);

\draw[fill=black] (0.606475,0.678096) circle (0.03);
\draw (0.606475,0.678096) node[anchor=south east,inner xsep=4pt,inner ysep=-1pt,scale=0.6] {\Large$5$};
\draw (0,0) -- (0.606475,0.678096);
\draw[path fading=north,fading transform={rotate=293.191}] (0.606475,0.678096) -- (0.836274,0.776546);
\draw[path fading=north,fading transform={rotate=343.191}] (0.606475,0.678096) -- (0.678769,0.917415);

\draw[fill=black] (0.937284,0.0445918) circle (0.03);
\draw (0.937284,0.0445918) node[anchor=south west,inner xsep=-2pt,inner ysep=7pt,scale=0.6] {\Large$6$};
\draw (0,0) -- (0.937284,0.0445918);
\draw[path fading=north,fading transform={rotate=272.724}] (0.937284,0.0445918) -- (1.187,0.0564722);
\draw[path fading=north,fading transform={rotate=232.724}] (0.937284,0.0445918) -- (1.13622,-0.106823);
\draw[path fading=north,fading transform={rotate=312.724}] (0.937284,0.0445918) -- (1.12094,0.214208);

\draw (0.627979,-0.655302) -- (0.937284,0.0445918);
\draw (0.937284,0.0445918) -- (0.606475,0.678096);

\end{tikzpicture}
\caption{Example of an isolated star around the vertex $i \in V$. Although the neighborhood consists of vertices $N(i) = \{1, \ldots, 6\}$, the largest isolated star is $S(i) = \{1, \ldots, 5\}$.}
\label{fig:isolated_star_example}
\end{figure}
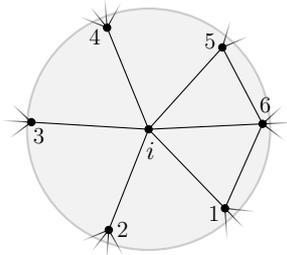

Intuitively, under $H_0$, the observed graph does not contain large isolated stars because of the underlying geometric structure. In fact, any isolated star under $H_0$ cannot be larger than the kissing number $\kissingnumber{d}$, which is the maximum number of non-overlapping spheres of the same radius that can be placed tangent to some central sphere in dimension $d$. To see this, note that our model is equivalent to the model where every vertex is the center of a sphere of radius $r / 2$, and two vertices are connected when their spheres touch or overlap. This means that, under $H_0$, it is impossible to observe an isolated star that is larger than the kissing number $\kissingnumber{d}$. For example, the kissing number for dimension $d = 2$ is $\kissingnumber{2} = 6$, so it is impossible to have more than six vertices in a given neighborhood without some of them being connected, see Figure~\ref{fig:isolated_star_example} for an example.

However, under the alternative hypothesis the observed graph can, and likely will, contain large isolated stars. In particular, a botnet vertex is quite likely to have an isolated star that is almost as large as its degree. Therefore, it will be likely to observe a few isolated stars that are larger than the kissing number $\kissingnumber{d}$. Hence, we can scan the graph and compute the size of the isolated star at every vertex. Then we reject $H_0$ when we see an isolated star that is larger than the kissing number.

\begin{definition}
\label{def:isolated_star_test}
Let $\kissingnumber{d}$ be the kissing number in dimension $d$. The \emph{isolated star test} rejects the null hypothesis for a given graph $G$ when $\max_{i \in V} |S(i)| > \kissingnumber{d}$.
\end{definition}

Checking whether there exists a vertex that has an isolated star that is larger than the kissing number can be done in $\bigO\bigl(\sum_{i \in V} d_i^{\smash{\kissingnumber{d}^2}}\bigr)$ time, with $d_i$ the degree of vertex $i \in V$. This scales polynomially in the number of vertices. However, in practice this is not feasible on large graphs, unless all vertices have quite small degree. Instead we can use a greedy algorithm to obtain lower bounds on the size of an isolated star, for example as described in \cite{Boppana1992}. Moreover, note that the kissing number $\kissingnumber{d}$ depends on the underlying dimension $d$, and the exact kissing number $\kissingnumber{d}$ is unknown for many dimensions. However, there exist good upper bounds which can be used instead. For dimensions $d \leq 24$, the best known upper bounds can be found in \cite{Mittelmann2010}, and for larger dimensions one could use the upper bound $\kissingnumber{d} \ll 1.3233^d$ \cite{Kabatiansky1978}.

Next we present the main result of this section, where we give conditions for the isolated star test to be asymptotically powerful. The proof of this result is postponed until Section~\ref{sec:proof_isolated_star_test_powerful}.
\begin{theorem}
\label{thm:isolated_star_test_powerful}
If $n p k \to \infty$ then the isolated star test from Definition~\ref{def:isolated_star_test} is asymptotically powerful, meaning that it has a risk converging to zero.
\end{theorem}

%%%%%%%%%%%%%%%%%%%%%%%%%%%%%%%%%%%%%%%%%%%%%%%%%%%%%%%%%%%%%%%%%%%%%%%%%%%%%%%%
%%%%%%%%%%%%%%%%%%%%%%%%%%%%%%%%%%%%%%%%%%%%%%%%%%%%%%%%%%%%%%%%%%%%%%%%%%%%%%%%
\subsubsection{Average distance test}
\label{sec:distance_based_test}
In this section we define a test that can detect whether an observed graph contains a planted botnet based on the difference in graph distances under the null and alternative hypothesis. Here we require that $p$ is large enough to ensure that the graph is connected with high probability.

Given two connected vertices $i, j \in V$, let $\distG(i, j)$ be the graph distance between $i$ and $j$. That is, $\distG(i, j)$ is the length of the shortest path in the graph $G$ that connects $i$ to $j$. Also, we define the average graph distance as
\begin{equation}
\distG^{\textup{avg}}(G) \coloneqq \frac{\sum_{1 \leq i < j \leq n} \1{\{i \leftrightsquigarrow j\}} \distG(i, j)}{\sum_{1 \leq i < j \leq n} \1{\{i \leftrightsquigarrow j\}}} \,.
%  \binom{n}{2}^{-1} \!\!\! \sum_{1 \leq i < j \leq n} \distG(i, j) \,.
\end{equation}

Under the null hypothesis, the observed graph is a random geometric graph and therefore the average graph distance will be large. To see this, consider first the average Euclidean distance between two uniformly chosen points on the torus. This can be lower bounded by
\begin{align}
\E_0[\distT(X_1, X_2)]
  &= \int_{[0, 1]^d} \sqrt{{\textstyle\sum_{j = 1}^d} \min\bigl(|x_j - 1/2|, 1 - |x_j - 1/2|\bigr)^2} \; \dif x_1 \cdots \dif x_d\\
  &= \int_{[0, 1]^d} \sqrt{{\textstyle\sum_{j = 1}^d} |x_j - 1/2|^2} \; \dif x_1 \cdots \dif x_d\\
\label{eq:average_torus_distance_lower_bound}
  &\geq \int_{[0, 1]^d} \: \max_{1 \leq j \leq d} |x_j - 1/2| \; \dif x_1 \cdots \dif x_d
  = \frac{d}{2(d+1)} \,,
\end{align}
where the final step follows by symmetry and is simply the expectation of the maximum of $d$ independent uniform random variables on $[0, 1/2]$. Hence, two uniformly chosen vertices have an expected Euclidean distance of at least $d / (2d + 2)$ on the torus. Then, consider the following lower bound on the average graph distance, which holds with high probability
\begin{equation}
\label{eq:average_distance_random_geometric_graph_lower_bound}
\distG^{\textup{avg}}(G)
  \geq \binom{n}{2}^{\!-1} \!\!\! \sum_{1 \leq i < j \leq n} \! \frac{\distT(X_i, X_j)}{r} \,,
\end{equation}
because we assumed that the graph is connected with high probability and because every edge can only connect two vertices when they are within distance $r$, so $\distG(i,j) \geq \distT(X_i,X_j) / r$. Note that, the right-hand side of \eqref{eq:average_distance_random_geometric_graph_lower_bound} can be seen as a U-statistic. Therefore, using \cite[Theorem 5.2]{Hoeffding1948}, we obtain
\begin{equation}
\text{Var}_0\left({\textstyle \binom{n}{2}^{-1} \sum_{1 \leq i < j \leq n} \distT(X_i, X_j)}\right)
  \leq \frac{2}{n} \, \text{Var}_0\left(\distT(X_1, X_2)\right)
  \to 0\,.
\end{equation}
Hence, Chebyshev's inequality ensures that $\binom{n}{2}^{-1} \sum_{1 \leq i < j \leq n} \distT(X_i, X_j)$ is concentrated around $\E_0[\distT(X_1, X_2)]$ with probability tending to one. Therefore, using \eqref{eq:average_torus_distance_lower_bound} and \eqref{eq:average_distance_random_geometric_graph_lower_bound}, we obtain for any $\epsilon > 0$ the following with high probability lower bound  
\begin{equation}
\label{eq:average_distance_random_geometric_graph_lower_bound2}
\distG^{\textup{avg}}(G)
  \geq (1 - \epsilon) \, \frac{\E_0[\distT(X_1, X_2)]}{r}\\
  \geq (1 - \epsilon) \, \frac{d}{2(d+1)} \cdot \frac{1}{r} \,,
\end{equation}

As we show below, the average graph distance is significantly smaller under the alternative hypothesis. Therefore, we consider the following test based on the average graph distance in the observed graph:
\begin{definition}
\label{def:average_distance_test}
Fix $\epsilon > 0$. The \emph{average distance test} rejects the null hypothesis for a given graph $G$ when
\begin{equation}
\label{eq:average_distance_test_boundary}
\distG^{\textup{avg}}(G)
  < (1 - \epsilon) \frac{d}{2(d+1)} \cdot \frac{1}{r} \,.
\end{equation}
\end{definition}

This brings us to the main result of this section, which identifies when the average distance test is asymptotically powerful. We postpone the proof of this theorem to Section~\ref{sec:proof_average_distance_test_powerful}.
\begin{theorem}
\label{thm:average_distance_test_powerful}
If $n p k \to \infty$ and $p$ is large enough to ensure that the subgraph induced by all non-botnet vertices is connected with high probability, then the average distance test from Definition~\ref{def:average_distance_test} is asymptotically powerful.
\end{theorem}

%\begin{remark}
Note that the assumption of connectedness implies that $n p \geq \bigOmega(\log(n))$ \cite{Penrose2003}, and together with the fact that $k \geq 1$ this implies $n p k \to \infty$. We include the latter condition to be able to compare the theorem above to Theorem~\ref{thm:isolated_star_test_powerful}. The requirement of connectivity is only a technical assumption that we make to considerably simplify the proof. This leads us to conjecture that Theorem~\ref{thm:average_distance_test_powerful} also holds under the milder condition that $p$ is large enough to ensure the existence of a giant component. This is also supported by our numerical simulations.
%\end{remark}

%%%%%%%%%%%%%%%%%%%%%%%%%%%%%%%%%%%%%%%%%%%%%%%%%%%%%%%%%%%%%%%%%%%%%%%%%%%%%%%%
%%%%%%%%%%%%%%%%%%%%%%%%%%%%%%%%%%%%%%%%%%%%%%%%%%%%%%%%%%%%%%%%%%%%%%%%%%%%%%%%
\subsubsection{Unknown dimension and connection radius}
\label{sec:unknown_dimension}
Computing the isolated star test requires knowledge of the dimension $d$ of the embedding space, and the average distance test requires the knowledge of the dimension $d$ as well as the connection radius $r$. In this section we show how to estimate these parameters from the observed graph.

%\paragraph{Estimating the dimension.}
%The problem of estimating the dimension $d$ was first studied in \cite{Serra2017}, where the authors considered a general setting that assumes no prior knowledge about the geometry. In our case, we can use a simpler approach because the geometry is known to be a $d$-dimensional torus. Therefore, 
To estimate the dimension $d$ we use the clustering coefficient \cite{Dall2002}. This is defined as the probability that two random neighbors of a given vertex are themselves connected. Under the null hypothesis, the clustering coefficient can be computed analytically and the resulting quantity only depends on the dimension $d$. Using \cite[see (15)]{Hammersley1950}, for distinct $i,j,k \in V$, we obtain 
\begin{align}
\label{eq:clustering_coefficient_analytic}
C_d
  &= \P_0(j \leftrightarrow k \,|\, i \leftrightarrow j, i \leftrightarrow k)\\
%  &= I_{3/4}\left(\frac{d+1}{2}, \frac{1}{2}\right) + I_{1/4}\left(\frac{d+1}{2}, \frac{d+1}{2}\right)\\
  &= \P\left(\text{Beta}\left(\frac{d+1}{2}, \frac{1}{2}\right) \leq \frac{3}{4}\right) + \P\left(\text{Beta}\left(\frac{d+1}{2}, \frac{d+1}{2}\right) \leq \frac{1}{4}\right) \,,
\end{align}
where $\text{Beta}(\cdot, \cdot)$ denotes a random variable with a beta distribution. Moreover, for a given graph, the clustering coefficient can be estimated by
\begin{equation}
\label{eq:clustering_coefficient_estimate}
\widehat{C}_d
  = \frac{\sum_{1 \leq i, j, k \leq n} \1{\{i \leftrightarrow j, i \leftrightarrow k, j \leftrightarrow k\}}}{\sum_{1 \leq i, j, k \leq n} \1{\{i \leftrightarrow j, i \leftrightarrow k\}}} \,.
\end{equation}
To estimate the dimension $d$ we can estimate the clustering coefficient $\widehat{C}_d$ using \eqref{eq:clustering_coefficient_estimate} and then invert the relation in \eqref{eq:clustering_coefficient_analytic} to obtain an estimate for the dimension $\widehat{d}$. This method of estimating the dimension gives a consistent estimator, under the null as well as the alternative hypothesis. This is shown in the following lemma, which we prove in Section~\ref{sec:proof_dimension_estimator_consistent}.
\begin{lemma}
\label{lem:dimension_estimator_consistent}
Using the clustering coefficient to estimate the dimension $\widehat{d}$ is consistent under both the null and alternative hypothesis, in the sense that both $\widehat{d} \conv[\P_0] d$ and $\widehat{d} \conv[\P_B] d$.
\end{lemma}

%\paragraph{Estimating the connection radius.}
The average distance test also requires knowledge of the connection radius $r$. To estimate this, we use that the edge probability $p$ is given by
\begin{equation}
\label{eq:edge_probability_analytic}
p = \P_0(i \leftrightarrow j) = \frac{(\sqrt{\pi} r)^d}{\Gamma(d/2 + 1)} \,,
\end{equation}%
where $\Gamma(\cdot)$ denotes the Gamma function. For a given graph, the edge probability can be estimated by
\begin{equation}
\label{eq:edge_probability_estimate}
\widehat{p} = \binom{n}{2}^{\!-1} \!\! \sum_{1 \leq i < j \leq n} \! \1{\{i \leftrightarrow j\}} \,.
\end{equation}

To obtain an estimate of the connection radius $r$, we can estimate the edge probability using \eqref{eq:edge_probability_estimate} and then invert the relation in \eqref{eq:edge_probability_analytic}, using our estimate of $d$, to obtain an estimate for the connection radius $\widehat{r}$. This method gives a consistent estimator of $p = p_n$, as the next lemma shows.
\begin{lemma}
\label{lem:probability_estimator_consistent}
Using $\widehat{p}$ to estimate $p = p_n$ is consistent both under the null and alternative hypothesis, in the sense that both $\widehat{p} / p \conv[\P_0] 1$ and $\widehat{p} / p \conv[\P_B] 1$.
\end{lemma}
We postpone the proof of Lemma \ref{lem:probability_estimator_consistent} to Section~\ref{sec:proof_connection_probability_estimator_consistent}. Since the radius $r$ is given in terms of a continuous function of $p$ in \eqref{eq:edge_probability_analytic}, this also shows that $\widehat{r}/r \conv[\P_0] 1$ and $\widehat{r}/r \conv[\P_B] 1$ by the continuous mapping theorem. Therefore, our estimate for the connection radius $\widehat{r}$ is also consistent under the null and alternative hypotheses.

\subsection{Identifying the botnet}
\label{sec:identifying_the_botnet}
When a test rejects the null hypothesis, we would also like to identify the vertices that are part of the botnet. To this end, let $\est \subseteq V$ be an estimator of the vertices in the botnet. We assume that the size of the botnet $|B| = k$ is known and that $|\est| = k$. To measure the performance of our estimator we use the risk function
\begin{equation}
\label{eq:risk_function}
\estrisk(\est) \coloneqq \E_B\left[\frac{|\est \symdiff B|}{2\,|B|}\right] \,,
\end{equation}
where $\est \symdiff B = ((V \setminus \est) \cap B) \cup (\est \cap (V \setminus B))$ is the symmetric difference between an estimator $\est$ of the botnet and the true botnet $B$. The reason for the normalization in \eqref{eq:risk_function} is that $|\est \symdiff B|$ could be unbounded, while  $0 \leq |\est \symdiff B|/|B| \leq 2$.

We say that a method achieves \emph{exact recovery} when $\estrisk(\est) \to 0$, and \emph{partial recovery} when $\estrisk(\est) \to \alpha$ for $\alpha \in (0, 1)$. In other words, partial recovery corresponds to identifying a positive proportion of the botnet vertices while exact recovery corresponds to identifying the majority of the botnet vertices. Note that, partial recovery is most interesting when the botnet size diverges. To see this, consider partial recovery of a single botnet vertex $k = 1$, in this case $\estrisk(\est) = \P_B(\est \neq B)$. That is, the botnet vertex is identified correctly only a fraction of the time, and remains unidentified otherwise.

Intuitively, our procedure identifies a botnet vertex when that vertex has a large enough isolated star $|S(i)|$. However, in this case, non-botnet vertices could have an isolated star that is larger than the kissing number $\kissingnumber{d}$, because it is connected to one or more botnet vertices. Therefore, in order to control the number of false positives, we introduce a parameter $\xi_n > 0$ to artificially increase the threshold $\kissingnumber{d}$ that was used when detecting the presence of a botnet. This leads to the following definition of the isolated star estimator:
\begin{definition}
\label{def:isolated_star_estimator}
Let $\kissingnumber{d}$ be the kissing number in dimension $d$. The \emph{isolated star estimator} is
\begin{equation}
\label{eq:isolated_star_estimator}
\est \coloneqq \left\{i \in V \,:\, |S(i)| > \kissingnumber{d} + \xi_n \right\} \,,
\end{equation}
with $\xi_n$ given by
\begin{equation}
\label{eq:isolated_star_estimator_threshold_increase}
\xi_n \coloneqq (1 + \epsilon) \frac{\log(n / k)}{\lambertW_0\left(\log(n / k) / (k p \e)\right)} \,,
\end{equation}
where $\epsilon > 0$ is arbitrary, and $\lambertW_0(\cdot)$ denotes the Lambert-W function\footnote{The function $\lambertW_0(\cdot)$ denotes one of the branches of the Lambert-W function. This is the solution in $y \in [-1, \infty)$ of the equation $x = y \mspace{1mu} \e^y$, with $x \geq -1 / \e$. For a detailed overview of this function and its properties see \cite{Corless1996}.}.
\end{definition}
Comparing this estimator with the isolated star test from Section~\ref{sec:isolated_star_test}, we see that the detection threshold is increased by $\xi_n$. In fact, we have chosen $\xi_n$ to be slightly larger than the maximum number of botnet vertices that are likely to connect to any non-botnet vertex. In other words, the addition of $\xi_n$ ensures that the number of false positives remains vanishingly small.

The performance of our test depends crucially on the asymptotic behavior of the expected number of edges $n p k$ that are connected to any botnet vertex. We will concisely refer to these as \textit{botnet edges}. Intuitively, when $n p k$ grows slowly, the botnet edges do not influence the largest isolated star of a typical vertex and thus $\xi_n$ is a constant. On the other hand, when $n p k$ is large, the largest isolated star of a typical vertex grows with $n$ and consequently $\xi_n$ also increases with $n$.

More precisely, we show that when $n p k \leq n^{\beta}$ with $\beta \in (0,1)$ our method always achieves at least partial recovery. This corresponds to the most common situation where there is a small botnet in a sparse graph. In this case, $\xi_n$ can be shown to converge to a constant, and thus every vertex with an isolated star that is only slightly larger than the kissing number $\kissingnumber{d}$ is considered a botnet vertex. On the other hand, if $n p k$ grows linearly in $n$ or faster, then the typical size of the largest isolated star is significantly larger than the kissing number $\kissingnumber{d}$ and additional technical assumptions are required for our method to achieve at least partial recovery. We make the above considerations precise in the main result of this section, which is presented below. 
\begin{theorem}
\label{thm:isolated_star_estimator_powerful}
Suppose that one of the following conditions holds:
\begin{enumerate}[label=(\roman*)]
\item\label{thm:isolated_star_estimator_powerful_part1}
$n p k \leq n^{\beta}$ for some $\beta \in (0,1)$,

\item\label{thm:isolated_star_estimator_powerful_part2}
$n^{1-\smallO(1)} \leq n p k \leq \smallO(n\log(n / k))$, and $\log(n / k)^2/n \leq p \leq \log(n / k)^{-2}$,

\item\label{thm:isolated_star_estimator_powerful_part3}
$n p k \geq \bigOmega(n\log(n / k))$, and $p = \smallO(k^{-2/3})$.
\end{enumerate}
Then the isolated star estimator from Definition~\ref{def:isolated_star_estimator} has exact recovery if $n p \to \infty$, and partial recovery otherwise.
\end{theorem}

Note that, when taken together, conditions \ref{thm:isolated_star_estimator_powerful_part1}--\ref{thm:isolated_star_estimator_powerful_part3} describe all possible asymptotic behaviors of $n p k$, but additional technical assumptions are required when $n p k \geq n^{1-\smallO(1)}$. The proof of Theorem~\ref{thm:isolated_star_estimator_powerful} is given in Section~\ref{sec:proof_isolated_star_estimator_powerful}.

%%%%%%%%%%%%%%%%%%%%%%%%%%%%%%%%%%%%%%%%%%%%%%%%%%%%%%%%%%%%%%%%%%%%%%%%%%%%%%%%
%%%%%%%%%%%%%%%%%%%%%%%%%%%%%%%%%%%%%%%%%%%%%%%%%%%%%%%%%%%%%%%%%%%%%%%%%%%%%%%%
\section{Simulations}
\label{sec:simulations}
We have shown that the tests introduced in the previous sections are asymptotically powerful when $n p k \to \infty$. In this section, we study the finite sample performance of these tests using simulations in order to compare their efficiency in practice on relatively small graphs. As specified both our tests have type-1 error that is nearly zero, so they will almost always correctly identify a graph without a botnet. Therefore, the focus of these simulations is on the type-2 error, which indicates how often a planted botnet is detected when it is actually present.

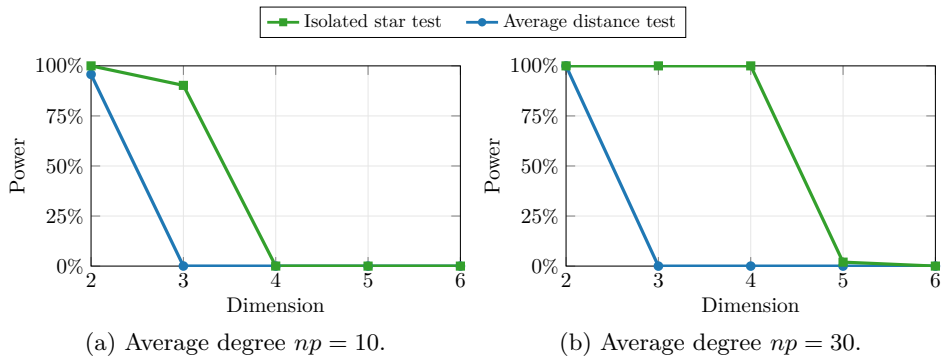
\begin{figure}[b!]
\definecolor{dark1}{RGB}{  31,120,180}
\definecolor{light1}{RGB}{149,188,209}
\definecolor{dark2}{RGB}{  51,160, 44}
\definecolor{light2}{RGB}{161,205,121}
\centering
\begin{tikzpicture}[scale=0.7,transform shape]
\begin{axis}[
hide axis,xmin=0,xmax=1,ymin=0,ymax=1,
legend cell align=left,legend columns=2]
% Isolated star test
\addlegendimage{color=dark2,line width=1.5pt,mark=square*,mark size=1.40pt}
\addlegendentry{Isolated star test$\quad$}
% Average distance test
\addlegendimage{color=dark1,line width=1.5pt,mark=*,mark size=1.50pt}
\addlegendentry{Average distance test}
\end{axis}
\end{tikzpicture}\\[4pt]%
\begin{subfigure}[t]{0.49\textwidth}
\scalebox{0.8}{\begin{tikzpicture}%[scale=0.8,transform shape]
\begin{axis}[
width=1.25\textwidth,
height=0.80\textwidth,
xlabel={Dimension},
ylabel={Power},
xmin=2,xmax=6,xtick={2,3,4,5,6},
ymin=0,ymax=100,ytick={0,25,50,75,100},
yticklabels={$0\text{\small\%}$,$25\text{\small\%}$,$50\text{\small\%}$,$75\text{\small\%}$,$100\text{\small\%}$},
xlabel style={yshift=2pt},
ylabel style={yshift=-2pt},
xticklabel style={/pgf/number format/fixed, /pgf/number format/fixed zerofill, /pgf/number format/precision=0},
grid=major,
grid style={draw=black!10!white}]
% Average distance test
\addplot[color=dark1,line width=1.5pt,mark=*,mark size=1.50pt] coordinates {
(2,  95.74) (3,  0.04) (4,  0.0) (5, 0.0) (6, 0.0)};
% Isolated star test
\addplot[color=dark2,line width=1.5pt,mark=square*,mark size=1.40pt] coordinates {
(2, 100.00) (3, 90.26) (4, 0.02) (5, 0.0) (6, 0.0)};
\end{axis}
\end{tikzpicture}}
\caption{Average degree $n p = 10$.}
\label{fig:simulation_results_normal_deg10}
\end{subfigure}%
\hspace{0.01\textwidth}%
\begin{subfigure}[t]{0.49\textwidth}
\scalebox{0.8}{\begin{tikzpicture}%[scale=0.8,transform shape]
\begin{axis}[
width=1.25\textwidth,
height=0.80\textwidth,
xlabel={Dimension},
ylabel={Power},
xmin=2,xmax=6,xtick={2,3,4,5,6},
ymin=0,ymax=100,ytick={0,25,50,75,100},
yticklabels={$0\text{\small\%}$,$25\text{\small\%}$,$50\text{\small\%}$,$75\text{\small\%}$,$100\text{\small\%}$},
xlabel style={yshift=2pt},
ylabel style={yshift=-2pt},
xticklabel style={/pgf/number format/fixed, /pgf/number format/fixed zerofill, /pgf/number format/precision=0},
grid=major,
grid style={draw=black!10!white}]
% Average distance test
\addplot[color=dark1,line width=1.5pt,mark=*,mark size=1.50pt] coordinates {
(2, 100.0) (3,   0.0) (4,   0.0) (5,  0.0) (6, 0.0)};
% Isolated star test
\addplot[color=dark2,line width=1.5pt,mark=square*,mark size=1.40pt] coordinates {
(2, 100.0) (3, 100.0) (4, 100.0) (5, 1.98) (6, 0.0)};
\end{axis}
\end{tikzpicture}}
\caption{Average degree $n p = 30$.}
\label{fig:simulation_results_normal_deg30}
\end{subfigure}%
\hspace{0.01\textwidth}%
\caption{The power of the isolated star test and the average distance test as a function of the dimension $d$. The threshold for rejecting the null hypothesis is as described in Sections \ref{sec:isolated_star_test} or \ref{sec:distance_based_test}, using estimated model parameters as described in Section~\ref{sec:unknown_dimension}. The parameters are: graph size $n = 10000$, botnet size $k = 10$, and each simulation contains $5000$ samples.}
\label{fig:simulation_results_normal}
\end{figure}

% Results of the simulation using our (theoretical) thresholds.
For our first simulation study we estimate the graph parameters with the consistent estimators described in Section~\ref{sec:unknown_dimension} and use these to compute the thresholds for rejecting the null hypothesis as explained in Sections \ref{sec:isolated_star_test} and \ref{sec:distance_based_test}. Further, for the isolated star test we use the greedy algorithm described in \cite{Boppana1992} to approximate the isolated star size of a given vertex.
The results of this can be seen in Figure~\ref{fig:simulation_results_normal}. Here we can see that both the isolated star test and average distance test perform quite well, even on relatively small graphs, provided that the underlying dimension is small. Nevertheless, the isolated star test performs better than the average distance test, especially when $n p$ is large.

% Type-1 error because of mis-estimation of graph parameters.
Note that using the estimated model parameters as described in Section~\ref{sec:unknown_dimension} instead of the true values could introduce some errors, which in turn could lead to our tests being incorrectly calibrated and result in a type-1 error that is too large. To investigate this issue we repeated the simulation with no botnet (i.e., $k = 0$). Both our tests were always correct and did not reject the null hypothesis in any of the trials. Furthermore, the dimension was correctly estimated in all cases. So it would have made no difference if we used the true dimension $d$ instead of the estimated value $\widehat{d}$. We did see some estimation errors for the dimension, but these were only present when the underlying dimension $d$ was larger than $10$.
%In that case, we observed some estimation errors with the dimension typically being underestimated, resulting in an estimated dimension $\widehat{d}$ that is slightly lower than the true value.
Moreover, for the average distance test we also need to estimate the connection radius. The errors introduced by using the estimator $\widehat{r}$ compared to the true value $r$ were minimal, and using $r$ instead of $\widehat{r}$ yields essentially the same performance as in Figure~\ref{fig:simulation_results_normal}.

% Explain that our (theoretical) thresholds are very conservative, try monte carlo (using true model parameters).
The results in Figure~\ref{fig:simulation_results_normal} show that both the isolated star test and average distance test can perform well even on relatively small graphs. However, we see that their performance quickly deteriorates as the dimension increases. This happens because the rejection thresholds as described in Sections \ref{sec:isolated_star_test} and \ref{sec:distance_based_test} are much too conservative. 

\begin{figure}[H]
\definecolor{dark1}{RGB}{  31,120,180}
\definecolor{light1}{RGB}{149,188,209}
\definecolor{dark2}{RGB}{  51,160, 44}
\definecolor{light2}{RGB}{161,205,121}
\centering
\begin{tikzpicture}[scale=0.7,transform shape]
\begin{axis}[
hide axis,xmin=0,xmax=1,ymin=0,ymax=1,
legend cell align=left,legend columns=2]
% Isolated star test
\addlegendimage{color=dark2,line width=1.5pt,mark=square*,mark size=1.40pt}
\addlegendentry{Isolated star test ($5\text{\small\%}$-{\small{}MC})$\quad$}
% Average distance test
\addlegendimage{color=dark1,line width=1.5pt,mark=*,mark size=1.50pt}
\addlegendentry{Average distance test ($5\text{\small\%}$-{\small{}MC})}
% Isolated star test 2
\addlegendimage{color=light2,line width=1.5pt,mark=triangle*,mark size=1.65pt}
\addlegendentry{Isolated star test ($0.1\text{\small\%}$-{\small{}MC})$\quad$}
% Average distance test 2
\addlegendimage{color=light1,line width=1.5pt,mark=pentagon*,mark size=1.75pt}
\addlegendentry{Average distance test ($0.1\text{\small\%}$-{\small{}MC})}
\end{axis}
\end{tikzpicture}\\[4pt]%
\begin{subfigure}[t]{0.49\textwidth}
\scalebox{0.8}{\begin{tikzpicture}%[scale=0.8,transform shape]
\begin{axis}[
width=1.25\textwidth,
height=0.80\textwidth,
xlabel={Dimension},
ylabel={Power},
xmin=2,xmax=18,xtick={2,6,10,14,18},
ymin=0,ymax=100,ytick={0,25,50,75,100},
yticklabels={$0\text{\small\%}$,$25\text{\small\%}$,$50\text{\small\%}$,$75\text{\small\%}$,$100\text{\small\%}$},
xlabel style={yshift=2pt},
ylabel style={yshift=-2pt},
xticklabel style={/pgf/number format/fixed, /pgf/number format/fixed zerofill, /pgf/number format/precision=0},
grid=major,
grid style={draw=black!10!white}]
% Extra info box referencing the histogram figure.
\draw[-,black!30!white,line width=0.75pt](axis cs:7,60.98) to[out=60,in=180] (axis cs:9,89.50);
\draw[-{Latex[length=6pt,width=4pt]},black!30!white,line width=0.75pt](axis cs:7,89.50) to (axis cs:11,89.50);
\node[anchor=west,rectangle,rounded corners=2pt,line width=0.75pt,inner ysep=4pt,inner xsep=4pt,outer xsep=3pt,fill=white,draw=black!30!white] at (axis cs:11,89.50) {\scriptsize\color{black!60!white} See also Figure~\mbox{\ref{fig:histogram_empirical_distributions}}.};
% Average distance test  (0.1%)
\addplot[color=light1,line width=1.5pt,mark=pentagon*,mark size=1.75pt] coordinates {
(2, 100.0) (3, 100.0) (4, 100.0) (5, 100.0) (6, 82.10) (7, 17.04) (8,  3.66) (9,  1.38) (10,  0.54) (11,  0.46) (12, 0.24) (13, 0.34) (14, 0.18) (15, 0.18) (16, 0.08) (17, 0.18) (18, 0.16) (19, 0.08) (20, 0.04) (21, 0.04) (22, 0.18)};
% Average distance test  (5.0%)
\addplot[color=dark1,line width=1.5pt,mark=*,mark size=1.50pt] coordinates {
(2, 100.0) (3, 100.0) (4, 100.0) (5, 100.0) (6, 97.68) (7, 60.98) (8, 27.16) (9, 16.20) (10, 11.74) (11,  9.34) (12, 7.86) (13, 7.30) (14, 6.84) (15, 6.58) (16, 6.64) (17, 5.86) (18, 6.70) (19, 6.26) (20, 5.90) (21, 5.46) (22, 5.62)};
% Isolated star test (0.1%)
\addplot[color=light2,line width=1.5pt,mark=triangle*,mark size=1.65pt] coordinates {
(2, 100.0) (3, 100.0) (4, 99.98) (5, 99.46) (6, 89.04) (7, 75.04) (8, 56.42) (9, 22.04) (10, 12.56) (11,  6.38) (12, 6.40) (13, 1.44) (14, 1.92) (15, 0.78) (16, 0.24) (17, 0.28) (18, 0.12) (19, 0.18) (20, 0.06) (21, 0.20) (22, 0.12)};
% Isolated star test (5.0%)
\addplot[color=dark2,line width=1.5pt,mark=square*,mark size=1.40pt] coordinates {
(2, 100.0) (3, 100.0) (4, 100.0) (5, 99.92) (6, 97.14) (7, 89.50) (8, 76.08) (9, 56.52) (10, 40.38) (11, 24.00) (12, 14.06) (13, 7.28) (14, 10.32) (15, 4.92) (16, 6.38) (17, 2.48) (18, 3.08) (19, 3.92) (20, 4.18) (21, 4.80) (22, 5.04)};
\end{axis}
\end{tikzpicture}}
\caption{Average degree $n p = 10$.}
\label{fig:simulation_results_mc_deg10}
\end{subfigure}%
\hspace{0.01\textwidth}%
\begin{subfigure}[t]{0.49\textwidth}
\scalebox{0.8}{\begin{tikzpicture}%[scale=0.8,transform shape]
\begin{axis}[
width=1.25\textwidth,
height=0.80\textwidth,
xlabel={Dimension},
ylabel={Power},
xmin=2,xmax=18,xtick={2,6,10,14,18},
ymin=0,ymax=100,ytick={0,25,50,75,100},
yticklabels={$0\text{\small\%}$,$25\text{\small\%}$,$50\text{\small\%}$,$75\text{\small\%}$,$100\text{\small\%}$},
xlabel style={yshift=2pt},
ylabel style={yshift=-2pt},
xticklabel style={/pgf/number format/fixed, /pgf/number format/fixed zerofill, /pgf/number format/precision=0},
grid=major,
grid style={draw=black!10!white}]
% Average distance test  (0.1%)
\addplot[color=light1,line width=1.5pt,mark=pentagon*,mark size=1.75pt] coordinates {
(2, 100.0) (3, 100.0) (4, 100.0) (5, 100.0) (6, 100.0) (7, 98.14) (8, 43.60) (9, 4.32) (10, 3.76) (11, 2.12) (12, 1.80) (13, 1.26) (14, 1.06) (15, 0.30) (16, 0.40) (17, 0.14) (18, 0.50) (19, 0.38) (20, 0.18) (21, 0.24) (22, 0.62)};
% Average distance test  (5.0%)
\addplot[color=dark1,line width=1.5pt,mark=*,mark size=1.50pt] coordinates {
(2, 100.0) (3, 100.0) (4, 100.0) (5, 100.0) (6, 100.0) (7, 100.0) (8, 85.88) (9, 51.80) (10, 29.80) (11, 23.54) (12, 21.18) (13, 16.06) (14, 13.24) (15, 13.20) (16, 12.38) (17, 11.36) (18, 11.10) (19, 10.30) (20, 9.48) (21, 10.56) (22, 9.92)};
% Isolated star test (0.1%)
\addplot[color=light2,line width=1.5pt,mark=triangle*,mark size=1.65pt] coordinates {
(2, 100.0) (3, 100.0) (4, 100.0) (5, 100.0) (6, 100.0) (7, 100.0) (8, 100.0) (9, 100.0) (10, 100.0) (11, 100.0) (12, 100.0) (13, 99.96) (14, 99.56) (15, 98.72) (16, 96.96) (17, 96.76) (18, 85.68) (19, 86.90) (20, 86.92) (21, 87.04) (22, 86.82)};
% Isolated star test (5.0%)
\addplot[color=dark2,line width=1.5pt,mark=square*,mark size=1.40pt] coordinates {
(2, 100.0) (3, 100.0) (4, 100.0) (5, 100.0) (6, 100.0) (7, 100.0) (8, 100.0) (9, 100.0) (10, 100.0) (11, 100.0) (12, 100.0) (13, 100.0) (14, 99.98) (15, 99.90) (16, 99.92) (17, 99.72) (18, 98.86) (19, 99.08) (20, 96.34) (21, 97.22) (22, 97.18)};
\end{axis}
\end{tikzpicture}}
\caption{Average degree $n p = 30$.}
\label{fig:simulation_results_mc_deg30}
\end{subfigure}%
\hspace{0.01\textwidth}%
\caption{The power of the isolated star test and the average distance test. The threshold for rejecting the null hypothesis is obtained by Monte Carlo calibration that ensures respectively $\alpha = 5\text{\small\%}$ and $\alpha = 0.1\text{\small\%}$ type-1 error, assuming that the dimension $d$ and connection radius $r$ are known. The parameters are: graph size $n = 10000$, botnet size $k = 10$, and each simulation contains $5000$ samples.}
\label{fig:simulation_results_mc}
\end{figure}
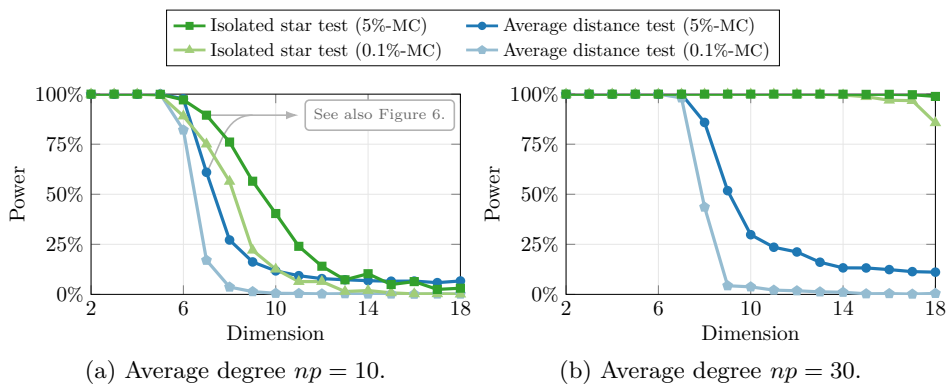

To better understand the properties of our two test statistics we conduct another simulation study, this time with clairvoyant knowledge of the dimension $d$ and connection radius $r$, which allows us to correctly calibrate these tests using a simple Monte Carlo method. That is, we sample $5000$ graphs from the null model (i.e., $k = 0$) and use these to compute the empirical distributions of either $\max_{i \in V} |S(i)|$ (for the isolated star test) and $\distG^{\textup{avg}}(G)$ (for the average distance test). We then take an appropriate quantile of these empirical distributions to obtain the rejection thresholds for a given significance level. The results of this can be seen in Figure~\ref{fig:simulation_results_mc}. This shows that the isolated star test outperforms the average distance test in most cases, especially when the dimension $d$ or the average degree $n p$ is large.

We note that the Monte Carlo method described above can also be applied when the dimension $d$ or the connection radius $r$ are unknown, but then using the estimated parameter as described in Section~\ref{sec:unknown_dimension}. However, the problem with this approach is that errors in the parameter estimation could lead to an incorrectly calibrated test, with a type-1 error that is possibly larger than the prescribed $\alpha$.

% Thorough explanation/comparison of our tests and results.
In Figure~\ref{fig:simulation_results_mc} we can see that the isolated star test has good performance when the dimension $d$ is small and the average degree $n p$ large. The reason for this is that the isolated star test rejects the null hypothesis when the graph contains an isolated star that is larger than a certain rejection threshold (i.e., the kissing number $\kissingnumber{d}$ in Figure~\ref{fig:simulation_results_normal}, or the threshold found by Monte Carlo calibration in Figure~\ref{fig:simulation_results_mc}). This rejection threshold is lower when the dimension $d$ is small, and the graph is more likely to contain a large isolated star when the average degree $n p$ is large. Hence, we see the best performance when the dimension $d$ is small and the average degree $n p$ large.

The performance of the average distance test is also related to the dimension $d$ and average degree $n p$ of the graph. To understand this, note that the botnet vertices can create shortcuts between vertices that are far away in the embedding space. When the average degree $n p$ is large, there is a higher probability that more shortcuts are created, which in turn decreases the average graph distance. On the other hand, as the dimension $d$ increases the average graph distance among the non-botnet vertices decreases, so the shortcuts created by any potential botnet vertices have a less pronounced effect. Thus, here we also see the best performance when the dimension $d$ is small and the average degree $n p$ large.

Finally, another reason why both tests have worse performance when the dimension $d$ increases is because the effect of the underlying geometry disappears when $d \to \infty$, as was shown in \cite{Bubeck2016}. Hence the difference between the null and alternative hypothesis is more pronounced when the dimension $d$ is small.

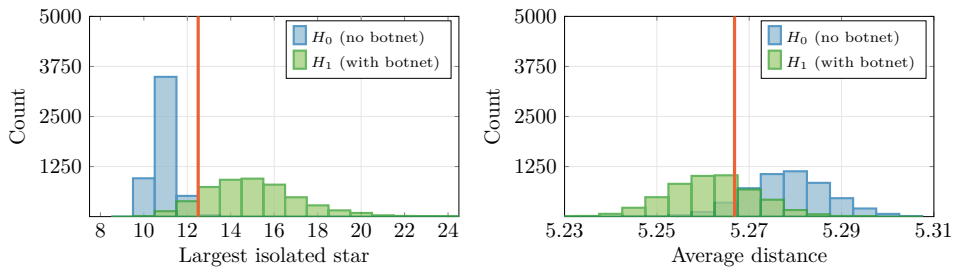
\begin{figure}[H]
\definecolor{dark1}{RGB}{  31,120,180}
\definecolor{light1}{RGB}{149,188,209}
\definecolor{dark2}{RGB}{  51,160, 44}
\definecolor{light2}{RGB}{161,205,121}
\definecolor{threshold1}{RGB}{235,98,53}
\centering
\begin{subfigure}[t]{0.49\textwidth}
\scalebox{0.8}{\begin{tikzpicture}%[scale=0.8,transform shape]
\begin{axis}[
ybar,
bar width=1,
bar shift=0pt,
ymax=1,ymin=0,
width=1.25\textwidth,
height=0.80\textwidth,
xlabel={Largest isolated star},
ylabel={Count},
xmin=7.5,xmax=24.5,xtick={8,10,14,18,22,24},extra x ticks={12,16,20},
ymin=0,ymax=5000,ytick={1250,2500,3750,5000},
xlabel style={yshift=2pt},
ylabel style={yshift=0pt},
xticklabel style={/pgf/number format/fixed, /pgf/number format/fixed zerofill, /pgf/number format/precision=0, /pgf/number format/1000 sep={}},
yticklabel style={/pgf/number format/fixed, /pgf/number format/fixed zerofill, /pgf/number format/precision=0, /pgf/number format/1000 sep={}},
%axis on top=true,
xtick align=inside,
major tick length=2pt,
grid=none,
legend cell align=left,
legend image code/.code={\draw [#1] (0cm,-0.1cm) rectangle (0.2cm,0.2cm); },
ymajorgrids,
grid style={draw=black!10!white},
extra tick style={
    grid=major,
    grid style={draw=black!10!white}
}]
%\node[text width=3cm,align=center,scale=0.9] at (axis cs: 9.1,3890) {\color{black!60!white}\small{}$H_0$ (no botnet)\\\raisebox{5pt}{\footnotesize{}$99.28${\scriptsize\%} correct}};
%\node[text width=3cm,align=center,scale=0.9] at (axis cs: 16.3,1348) {\color{black!60!white}\small{}$H_1$ (with botnet)\\\raisebox{5pt}{\footnotesize{}$89.50${\scriptsize\%} correct}};
\addplot[color=dark1,fill=light1,line width=1.0pt,opacity=0.75] coordinates {(9, 0) (10, 957) (11, 3490) (12, 517) (13, 36) (14, 0) (15, 0) (16, 0) (17, 0) (18, 0) (19, 0) (20, 0) (21, 0) (22, 0) (23, 0) (24, 0) (25, 0)};
\addplot[color=dark2,fill=light2,line width=1.0pt,opacity=0.75] coordinates {(9, 0) (10, 6) (11, 135) (12, 384) (13, 739) (14, 920) (15, 948) (16, 798) (17, 483) (18, 281) (19, 156) (20, 92) (21, 32) (22, 17) (23, 8) (24, 1) (25, 0)};
\draw[-,threshold1,line width=1.5pt](axis cs:12.5,0)--(axis cs:12.5,5000);
\legend{{\scriptsize{}$H_0$ (no botnet)}, {\scriptsize{}$H_1$ (with botnet)}}
\end{axis}
\end{tikzpicture}}
\end{subfigure}%
\hspace{0.01\textwidth}%
\begin{subfigure}[t]{0.49\textwidth}
\scalebox{0.8}{\begin{tikzpicture}%[scale=0.8,transform shape]
\begin{axis}[
ybar,
bar width=0.005,
bar shift=0pt,
ymax=1,ymin=0,
width=1.25\textwidth,
height=0.80\textwidth,
xlabel={Average distance},
ylabel={Count},
xmin=5.23,xmax=5.31,xtick={5.23,5.25,5.27,5.29,5.31},
ymin=0,ymax=5000,ytick={1250,2500,3750,5000},
xlabel style={yshift=2pt},
ylabel style={yshift=0pt},
xticklabel style={/pgf/number format/fixed, /pgf/number format/fixed zerofill, /pgf/number format/precision=2, /pgf/number format/1000 sep={}},
yticklabel style={/pgf/number format/fixed, /pgf/number format/fixed zerofill, /pgf/number format/precision=0, /pgf/number format/1000 sep={}},
%axis on top=true,
xtick align=inside,
major tick length=2pt,
grid=none,
legend cell align=left,
legend image code/.code={\draw [#1] (0cm,-0.1cm) rectangle (0.2cm,0.2cm); },
grid=major,
grid style={draw=black!10!white}]
%\node[text width=3cm,align=center,scale=0.9] at (axis cs: 5.2815,1532) {\color{black!60!white}\small{}$H_0$ (no botnet)\\\raisebox{5pt}{\footnotesize{}$95.00${\scriptsize\%} correct}};
%\node[text width=3cm,align=center,scale=0.9] at (axis cs: 5.2505,1431) {\color{black!60!white}\small{}$H_1$ (with botnet)\\\raisebox{5pt}{\footnotesize{}$60.98${\scriptsize\%} correct}};
\addplot[color=dark1,fill=light1,line width=1.0pt,opacity=0.75] coordinates {(5.23, 0) (5.235, 0) (5.24, 0) (5.245, 1) (5.25, 6) (5.255, 37) (5.26, 116) (5.265, 348) (5.27, 709) (5.275, 1062) (5.28, 1132) (5.285, 843) (5.29, 460) (5.295, 203) (5.3, 68) (5.305, 14)};
\addplot[color=dark2,fill=light2,line width=1.0pt,opacity=0.75] coordinates {(5.23, 1) (5.235, 17) (5.24, 70) (5.245, 221) (5.25, 485) (5.255, 818) (5.26, 1016) (5.265, 1031) (5.27, 678) (5.275, 422) (5.28, 164) (5.285, 59) (5.29, 16) (5.295, 2) (5.3, 0) (5.305, 0)};
\draw[-,threshold1,line width=1.5pt](axis cs:5.26682,0)--(axis cs:5.26682,5000);
\legend{{\scriptsize{}$H_0$ (no botnet)}, {\scriptsize{}$H_1$ (with botnet)}}
\end{axis}
\end{tikzpicture}}
\end{subfigure}%
\hspace{0.01\textwidth}%
\caption{Histograms comparing the empirical distributions of the largest isolated star $\max_{i \in V} |S(i)|$ and the average distance $\distG^{\textup{avg}}(G)$ statistics under the null and alternative hypothesis. The threshold for rejecting the null hypothesis at the $\alpha = 0.05$ significance level is shown in red. The parameters are: graph size $n = 10000$, botnet size $k = 10$, average degree $n p = 10$, dimension $d = 7$, and each histogram contains $5000$ samples.}
\label{fig:histogram_empirical_distributions}
\end{figure}

%%%%%%%%%%%%%%%%%%%%%%%%%%%%%%%%%%%%%%%%%%%%%%%%%%%%%%%%%%%%%%%%%%%%%%%%%%%%%%%%
%%%%%%%%%%%%%%%%%%%%%%%%%%%%%%%%%%%%%%%%%%%%%%%%%%%%%%%%%%%%%%%%%%%%%%%%%%%%%%%%
\section{Discussion}
\label{sec:discussion}
In this section we remark on our results and discuss some possible directions for future work.

\paragraph{Different null hypothesis.} 
Our results show that it is possible to detect an arbitrarily small planted botnet, provided that $n p k \to \infty$. However, these results hinge on the underlying geometric structure of the model. Many other network models have been developed that are based on a different geometry than the one assumed by our model \cite{Boguna2010,Krioukov2010,Barthelemy2011,Deijfen2011,Bringmann2019}. Therefore, it would be interesting to see what the effect of the underlying geometry is, and to what extent our results can be extended to models that have a different underlying geometric structure.

Our tests and analytical approach is fairly robust against minor changes in the underlying geometry. For instance, our results remain true when the embedding space is a slightly deformed torus or sphere, or the points are distributed in the embedding space in a slightly non-uniform way. However, when the changes in geometry are more drastic we expect the nature of the results to change. In particular, when the geometry causes the resulting graph to become a small world we expect the average distance test to fail, and when the geometry causes considerable inhomogeneity in vertex degrees we expect the isolated star test to fail.

\paragraph{Smaller isolated stars for higher power.}
The isolated star test rejects the null hypothesis when the largest observed isolated star is bigger than the kissing number $\kissingnumber{d}$, which automatically ensures that the type-1 error is zero. However, for dimensions $d > 2$, the typical largest isolated star in a random geometric graph is much smaller than the kissing number $\kissingnumber{d}$. For example, numerical simulations suggest that in dimension $d = 4$, the size of the typical isolated star is smaller than $10$, whereas the kissing number is $\kissingnumber{4} = 24$ \cite{Musin2003,Mittelmann2010}. This suggests that, depending on the significance level, one might use a much smaller threshold value, which would greatly increase the power of the test.

One possible way to achieve this is to calibrate the test using a Monte Carlo approach, as we did in Section~\ref{sec:simulations}. However, this is a computationally expensive approach which could be avoided with better knowledge of the behavior of isolated star sizes in higher dimensions.

\paragraph{Diverging dimension.}
From a theoretical perspective it would be interesting to know whether our results can be extended to the setting where the dimension $d$ is diverging together with the graph size $n$, similar to the problem considered in \cite{Bubeck2016}. For the isolated star test, we can use the following bound on the kissing number $\kissingnumber{d} \ll 1.3233^d$ \cite{Kabatiansky1978}. In this case, the same arguments as in the proof of Theorem~\ref{thm:isolated_star_test_powerful} suggest that the isolated star test is asymptotically powerful when $1 \ll n p \ll n^{1/3}$ and
\begin{equation}
\label{eq:growing_dimension_isolated_star_test_powerful}
d \leq \frac{\log(n p)}{\log(1.3233)} \,.
\end{equation}
%Hence, the isolated star test is asymptotically powerful even when the dimension diverges slowly, provided that the average degree $n p$ diverges as well.
However, a better understanding of the distribution of isolated stars in graphs with large underlying dimension could significantly improve this result and possibly show that the isolated star test can still be applied even when the dimension grows much faster than \eqref{eq:growing_dimension_isolated_star_test_powerful}.

\paragraph{Estimating the botnet size.}
In Section~\ref{sec:identifying_the_botnet} we show that, under some technical conditions, it is possible to asymptotically identify all botnet vertices provided $n p \to \infty$, and that a part of the botnet can be recovered when $n p = \bigO(1)$. It could be an interesting possibility for future research to see whether it is possible to estimate the botnet size $|B|$. In the setting where we have exact recovery (i.e., $n p \to \infty$) this is of course trivial, but it would be very interesting to see how well that botnet size $|B|$ can be estimated when $n p = \bigO(1)$.

%%%%%%%%%%%%%%%%%%%%%%%%%%%%%%%%%%%%%%%%%%%%%%%%%%%%%%%%%%%%%%%%%%%%%%%%%%%%%%%%
%%%%%%%%%%%%%%%%%%%%%%%%%%%%%%%%%%%%%%%%%%%%%%%%%%%%%%%%%%%%%%%%%%%%%%%%%%%%%%%%
\section{Proofs}
\label{sec:proofs}
This section is devoted to the proofs of the results stated in Sections \ref{sec:detecting_a_botnet} and \ref{sec:identifying_the_botnet}.

%%%%%%%%%%%%%%%%%%%%%%%%%%%%%%%%%%%%%%%%%%%%%%%%%%%%%%%%%%%%%%%%%%%%%%%%%%%%%%%%
%%%%%%%%%%%%%%%%%%%%%%%%%%%%%%%%%%%%%%%%%%%%%%%%%%%%%%%%%%%%%%%%%%%%%%%%%%%%%%%%
\subsection{Proof of Theorem~\ref{thm:isolated_star_test_powerful}: Isolated star test is powerful}
\label{sec:proof_isolated_star_test_powerful}
\begin{proof}[\unskip\nopunct]
As explained in Section~\ref{sec:isolated_star_test}, the isolated star test has zero type-1 error (i.e., it always correctly identifies a random geometric graph without a botnet). Therefore, to show that the isolated star test is asymptotically powerful, we must show that under the alternative hypothesis, the probability of having an isolated star larger than the kissing number $\kissingnumber{d}$ tends to one. This is done in two steps. First, let $\partdeg(i)$ be the non-botnet degree of a vertex $i \in V$. That is, $\partdeg(i)$ denotes the number of non-botnet neighbors of $i$. Then, we show that any botnet vertex $i \in B$, with $\partdeg(i) \geq \kissingnumber{d} + 1$, will form an isolated star of size $|S(i)| \geq \kissingnumber{d} + 1$ with high probability. Second, we show that, with high probability, there exists a botnet vertex that has arbitrarily large non-botnet degree.

Given a botnet vertex $i \in B$, define the event $\eventmd(i) \coloneqq \{\partdeg(i) \geq \kissingnumber{d} + 1\}$. Then, conditionally on the event $\eventmd(i)$, let $\{v_1, \ldots, v_{\kissingnumber{d}+1}\}$ be a subset of $\kissingnumber{d} + 1$ non-botnet neighbors of $i$. We reveal these vertices one at a time. For every vertex $v_j$ revealed this way, let $q_j$ be the probability that $v_j$ is not connected to any of the previously revealed vertices given that all these previously revealed vertices are themselves not connected. For $j \in [\kissingnumber{d}+1] = \{1, \ldots, \kissingnumber{d}+1\}$ we obtain
\begin{align}
q_j
  \coloneqq{}& \P_B\bigl(v_j \nleftrightarrow v_k \;\forall\, k \in [j-1]
    \,\big|\, \eventmd(i),\, v_k \nleftrightarrow v_l \;\forall\, k < l \in [j-1]\bigr)\notag\\
  ={}& 
  \P_B\bigl(\distT(X_{v_j}, X_{v_k}) > r \;\forall\, k \in [j-1]
      \,\big|\, \eventmd(i),\, v_k \nleftrightarrow v_l \;\forall\, k < l \in [j-1]\bigr)\notag\\
\label{eq:not_connected_previous_probability}
  \geq{}& 1 - (j-1) p \,,
\end{align}
where we note that, because $i \in B$ is a botnet vertex, conditioning on the event $\eventmd(i)$ does not affect the distribution of the vertex locations (i.e., these remain uniform random variables on the torus). Furthermore, observe that \eqref{eq:not_connected_previous_probability} becomes an equality precisely when the torus distance between every pair of previously revealed vertices is larger than $2 r$. Then, a lower bound on the probability that $i \in B$ forms an isolated star of size at least $\kissingnumber{d} + 1$ is given by
\begin{align}
\P_B\bigl(|S(i)| \geq \kissingnumber{d} + 1 \,\big|\, \eventmd(i)\bigr)
  &\geq \P_B\bigl(v_j \nleftrightarrow v_k \;\forall\, j < k \in [\kissingnumber{d}+1] \,\big|\, \eventmd(i)\bigr)\\
\label{eq:isolated_star_probability_lower_bound}
  &= \prod_{j=1}^{\kissingnumber{d}+1} q_j
  \geq (1-\kissingnumber{d} p)^{\kissingnumber{d}}
  \to 1 \,,
\end{align}
where the convergence to $1$ follows because $p \to 0$ and $\kissingnumber{d}$ is constant. Hence, any botnet vertex $i \in B$ with $\partdeg(i) \geq \kissingnumber{d} + 1$ will form an isolated star of size $|S(i)| \geq \kissingnumber{d} + 1$ with probability tending to one.

For the second part of the proof, we will show that there indeed exists a botnet vertex $i \in B$ with $\partdeg(i) \geq \kissingnumber{d} + 1$. First observe that for all $i \in B$ the non-botnet degrees $\partdeg(i)$ are independent random variables distributed as $\text{Bin}(n - k, p)$. Moreover, by the Stein-Chen method \cite{Chen1975,DenHollander2012}, it follows that 
\begin{equation}
\label{eq:stein_chen_partdeg_approximation}
\bigl\Vert{} \partdeg(i) - \text{Poi}((n-k) p) \bigr\Vert_{\text{\scriptsize{}TV}}
  \leq 2 p
  \to 0 \,,
\end{equation}
where $\smash{\Vert\mspace{-2mu}\cdot\mspace{-2mu}\Vert_{\text{\scriptsize{}TV}}}$ denotes the total variation norm. Now, because $n p k \to \infty$ and $k = \smallO(n)$ it follows that either $(n-k) p \to \infty$, or $(n-k) p = \bigTheta(1)$ and $k \to \infty$. When $(n-k) p \to \infty$ every botnet vertex will eventually have non-botnet degree larger than $\kissingnumber{d} + 1$ with high probability. On the other hand, if $(n-k) p = \bigTheta(1)$ then by \eqref{eq:stein_chen_partdeg_approximation} there is a positive probability that $\partdeg(i) \geq \kissingnumber{d} + 1$, independently for each botnet vertex $i \in B$, and since $k \to \infty$ there exists a botnet vertex with non-botnet degree larger than $\kissingnumber{d} + 1$ with high probability. Finally, combining this with \eqref{eq:isolated_star_probability_lower_bound} shows that the graph will contain an isolated star larger than $\kissingnumber{d} + 1$ with high probability.
\end{proof}

%%%%%%%%%%%%%%%%%%%%%%%%%%%%%%%%%%%%%%%%%%%%%%%%%%%%%%%%%%%%%%%%%%%%%%%%%%%%%%%%
%%%%%%%%%%%%%%%%%%%%%%%%%%%%%%%%%%%%%%%%%%%%%%%%%%%%%%%%%%%%%%%%%%%%%%%%%%%%%%%%
\subsection{Proof of Theorem~\ref{thm:average_distance_test_powerful}: Average distance test is powerful}
\label{sec:proof_average_distance_test_powerful}
\begin{proof}[\unskip\nopunct]
As given in \eqref{eq:average_distance_test_boundary}, under the null hypothesis we have the high probability lower bound
\begin{equation}
\distG^{\textup{avg}}(G)
  \geq (1 - \epsilon) \, \frac{d}{2(d+1)} \cdot \frac{1}{r} \,,
\end{equation}
Therefore, the average distance test has vanishing type-1 error (i.e., it will correctly identify a geometric random graph with no botnet with high probability). To show that this test is asymptotically powerful, we are left to show that the type-2 error also vanishes. This is done by showing that, under the alternative, there is a botnet vertex that creates a shortcut between most pairs of non-botnet vertices, as shown in Figure~\ref{fig:example_botnet_shortcut}. Using this, we show that, with high probability, the average graph distance is at most $\smallO(1) / r$, which is much smaller than the threshold in \eqref{eq:average_distance_test_boundary}.

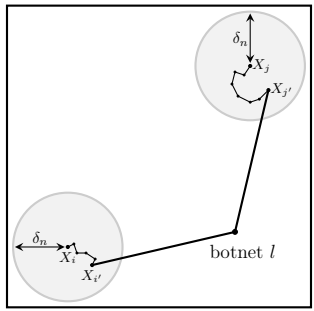
\begin{figure}[H]
\centering
\begin{tikzpicture}[scale=4.0]
% Frame
\draw[black,line width=0.8pt] (0.00,0.00) rectangle (1.00,1.00);
% Bottom left ball
\draw[white!80!black,line width=0.8pt,fill=white!95!black] (0.20,0.20) circle (0.18);
\draw[black,thin,stealth-stealth] (0.026,0.20) -- (0.190,0.20);
\draw (0.11,0.20) node[black,anchor=south,inner sep=1pt,scale=0.6] {$\delta_n$};

\draw[black] (0.20,0.20) -- (0.22,0.21) -- (0.23,0.18) -- (0.26,0.18) -- (0.29,0.16) -- (0.28,0.14);
\draw (0.20,0.20) node[anchor=north,inner xsep=3pt,inner ysep=3pt,scale=0.6] {\small{}$X_i$};
\draw (0.28,0.13) node[anchor=north,inner xsep=3pt,inner ysep=1pt,scale=0.6] {\small{}$X_{i'}$};
\draw[fill=black] (0.20,0.20) circle (0.005);
\draw[fill=black] (0.22,0.21) circle (0.003);
\draw[fill=black] (0.23,0.18) circle (0.003);
\draw[fill=black] (0.26,0.18) circle (0.003);
\draw[fill=black] (0.29,0.16) circle (0.003);
\draw[fill=black] (0.28,0.14) circle (0.005);

% Top right ball
\draw[white!80!black,line width=0.8pt,fill=white!95!black] (0.80,0.80) circle (0.18);
\draw[black,thin,stealth-stealth] (0.80,0.810) -- (0.80,0.974);
\draw (0.80,0.89) node[black,anchor=east,inner sep=0.75pt,scale=0.6] {$\delta_n$};

\draw[black] (0.80,0.80) -- (0.78,0.77) -- (0.75,0.78) -- (0.74,0.74) -- (0.76,0.70) -- (0.80,0.68) -- (0.83,0.69) -- (0.86,0.72);
\draw (0.80,0.80) node[anchor=west,inner xsep=1pt,inner ysep=0pt,scale=0.6] {\small{}$X_j$};
\draw (0.86,0.72) node[anchor=west,inner xsep=1pt,inner ysep=0pt,scale=0.6] {\small{}$X_{j'}$};
\draw[fill=black] (0.80,0.80) circle (0.005);
\draw[fill=black] (0.78,0.77) circle (0.003);
\draw[fill=black] (0.75,0.78) circle (0.003);
\draw[fill=black] (0.74,0.74) circle (0.003);
\draw[fill=black] (0.76,0.70) circle (0.003);
\draw[fill=black] (0.80,0.68) circle (0.003);
\draw[fill=black] (0.83,0.69) circle (0.003);
\draw[fill=black] (0.86,0.72) circle (0.005);
% Botnet vertex and its connections
%\draw[ultra thin,black] (0.28,0.14) .. controls (0.60,0.13) and (0.60,0.13) .. (0.75,0.25);
%\draw[ultra thin,black] (0.86,0.71) .. controls (0.87,0.40) and (0.87,0.40) .. (0.75,0.25);
%\draw[thick,dashed] (0.28,0.14) .. controls (0.60,0.13) and (0.60,0.13) .. (0.75,0.25);
%\draw[thick,dashed] (0.86,0.71) .. controls (0.87,0.40) and (0.87,0.40) .. (0.75,0.25);
\draw[thick,black] (0.28,0.14) -- (0.75,0.25);
\draw[thick,black] (0.86,0.72) -- (0.75,0.25);
\draw[fill=black] (0.75,0.25) circle (0.008);
\draw (0.75,0.25) node[anchor=north west,inner xsep=-16pt,inner ysep=8pt,align=left,scale=0.6] {\large{}botnet $l$};
\end{tikzpicture}
\caption{Example of botnet vertex $l \in B$ creating a shortcut between vertices $i, j \in V \setminus B$.}
\label{fig:example_botnet_shortcut}
\end{figure}

For a non-botnet vertex $i \in V \setminus B$, let $\ball(X_i; \delta_n)$ denote the ball of radius $\delta_n \coloneqq (V_d \log(n p))^{-1/d}$ around the location $X_i$, where $V_d \coloneqq \pi^{d/2} / \Gamma(d/2+1)$ denotes the volume of a $d$-dimensional unit ball. Also, let $\ballsubset_i \subseteq V \setminus B$ denote the non-botnet vertices with location in $\ball(X_i; \delta_n)$, that is
\begin{equation}
\ballsubset_i \coloneqq \left\{i' \in V \setminus B : X_{i'} \in \ball(X_i; \delta_n) \right\} \,.
\end{equation}
Note that, because $k = \smallO(n)$, we have
\begin{equation}
\E_B[|\ballsubset_i|]
  = \sum_{i' \in V \setminus B} \P_B(i' \in \ballsubset_i)
  = (n - k) V_d \delta_n^d
  = \frac{n - k}{\log(n p)}
  = (1 + \smallO(1)) \frac{n}{\log(n p)} \,.
\end{equation}
Therefore, using the relative Chernoff bound \cite[see (7)]{Hagerup1990} or \cite[Theorem~4.5]{Mitzenmacher2017}, for any $\xi > 0$, we obtain
\begin{align}
\P_B\left(|\ballsubset_i| \geq (1 - \xi) \frac{n}{\log(n p)}\right)
  &= 1 - \P_B\left(|\ballsubset_i| < (1 - \xi) \frac{n}{\log(n p)}\right)\\
  &\geq 1 - \P_B\bigl(|\ballsubset_i| < (1 - \xi / 2) \E[|\ballsubset_i|]\bigr)\\
  &\geq 1 - \exp\left(- \frac{\xi^2 \, n}{8 \log(n p)}\right)
  \to 1 \,.
\end{align}
Now, let $l \in B$ be an arbitrary botnet vertex, and consider the probability that there exists a vertex $i' \in \ballsubset_i$ that connects to the botnet vertex $l$. This gives
\begin{multline}
\P_B\bigl(\exists i' \in \ballsubset_i : i' \leftrightarrow l\bigr)\\
\begin{aligned}[b]
  &\geq \P_B\left(\exists i' \in \ballsubset_i : i' \leftrightarrow l \,\middle|\, |\ballsubset_i| \geq (1 - \xi)\frac{n}{\log(n p)}\right) \, \P_B\left(|\ballsubset_i| \geq (1 - \xi) \frac{n}{\log(n p)}\right) \notag\\
  &\geq (1 + \smallO(1)) \left(1 - (1 - p)^{(1 - \xi) n / \log(n p)}\right)\\
%\label{eq:botnet_shortcut_probability}
  &\geq (1 + \smallO(1)) \left(1 - \e^{-\xi n p / \log(n p)}\right)
  \to 1 \,,
\end{aligned}
\end{multline}
where the convergence to $1$ follows because $n p / \log(n p) \to \infty$. To continue, we use an existing result relating the torus distance and the graph distance \cite{Friedrich2013,Ellis2007,Muthukrishnan2005,Bradonjic2010,Diaz2016}. Translated to our notation, this result is as follows:
\begin{theorem*}[{see \cite[Theorem~3]{Friedrich2013} or \cite[Theorem~8]{Ellis2007}}]
There exists a constant $K$ independent of $n$ such that for any pair of vertices in the same connected component $i, j \in V$ with $\distT(X_i, X_j) \gg \frac{\log(n)}{n \, r^{d-1}}$ we obtain $\distG(i, j) \leq K \mspace{1mu} \distT(X_i, X_j) / r$ with high probability.
\end{theorem*}
Define the event $\eventcg \coloneqq \{G_{\scriptscriptstyle V \setminus B} \text{ is connected}\}$, where $G_{\scriptscriptstyle V \setminus B}$ denotes the subgraph induced by all non-botnet vertices. Note that $\P_B(\eventcg) \to 1$ by assumption. Then, given the event $\eventcg$, the result above guarantees that there exists a path of length at most $\bigO(\delta_n) / r$ between $i$ and every $i' \in \ballsubset_i$. Hence, for a given $i \in V \setminus B$,
\begin{multline}
\P_B\bigl(\distG(i, l) \leq 1 + \bigO(\delta_n) / r\bigr)\\
\begin{aligned}[b]
  &= \P_B\bigl(\eventcg \cap \bigl\{\distG(i, l) \leq 1 + \bigO(\delta_n) / r\bigr\}\bigr) - \smallO(1)\\
\label{eq:distance_to_botnet_vertex}
  &\geq \P_B\bigl(\eventcg \cap \bigl\{\exists i' \in \ballsubset_i : i' \leftrightarrow l, \distG(i, i') \leq \bigO(\delta_n) / r\bigr\}\bigr) - \smallO(1)
%  &\geq 1 - \bigl(1 - \P_B\bigl(\eventcg)\bigr) - \bigl(1 - \P_B\bigl(\exists i' \in \ballsubset_i : i' \leftrightarrow l, \distG(i, i') \leq \bigO(\delta_n) / r\bigr)\bigr) - \smallO(1)
  \to 1 \,.
\end{aligned}
\end{multline}
Then, by definition of $\delta_n$ and applying \eqref{eq:distance_to_botnet_vertex} twice, we obtain for an arbitrary pair of non-botnet vertices $i, j \in V \setminus B$ and botnet vertex $l \in B$,
\begin{multline}
\label{eq:nonbotnet_distance_bound}
\P_B\bigl(\distG(i, j) \leq \smallO(1) / r\bigr)\\
\begin{aligned}[b]
  &\geq \P_B\bigl(\distG(i, j) \leq 2 + 2 \, \bigO(\delta_n) / r\bigr)\\
  &\geq \P_B\bigl(\distG(i, l) \leq 1 + \bigO(\delta_n) / r,\, \distG(j, l) \leq 1 + \bigO(\delta_n) / r\bigr)
  \to 1 \,.
\end{aligned}
\end{multline}
By observing that every botnet vertex connects to several non-botnet vertices with high probability (as explained at the end of the proof of Theorem~\ref{thm:isolated_star_test_powerful}), the above can be strengthened to also include the botnet vertices, and show that the distance between any given pair of vertices is at most $\smallO(1) / r$ with high probability. This brings us to the central result of this proof, namely that for an arbitrary pair $i, j \in V$ it follows that
\begin{equation}
\label{eq:distance_bound}
\P_B\bigl(\distG(i, j) \leq \smallO(1) / r\bigr) \to 1 \,,
\end{equation}

We continue by showing that the diameter of the graph $G$ is at most $\bigO(1) / r$ with high probability. To this end, we first consider the diameter of $G_{\scriptscriptstyle V \setminus B}$, this gives
\begin{multline}
\P_B\Bigl(\max_{i, j \in V \setminus B} \distG(i, j) \leq \bigO(1) / r\Bigr)\\[-2pt]
\begin{aligned}[b]
  &= \P_B\bigl(\diam(G_{\scriptscriptstyle V \setminus B}) \leq \bigO(1) / r\bigr)\\
\label{eq:nonbotnet_diameter_bound}
  &= \P_B\bigl(\eventcg \cap \bigl\{\diam(G_{\scriptscriptstyle V \setminus B}) \leq \bigO(1) / r\bigr\}\bigr) - \smallO(1)
  \to 1 \,,
\end{aligned}
\end{multline}
where the convergence to $1$ follows from the theorem stated above (see also \cite[Corollary~6]{Friedrich2013}). Similarly to what we did above, this can be extended to the diameter of $G$ by showing that every botnet vertex connects to at least one non-botnet vertex. Let $l \in B$ denote an arbitrary botnet vertex, then
\begin{align}
\P_B\Bigl(\min_{i \in B} \, \partdeg(i) \geq 1\Bigr)
  &= 1 - \bigl(\P_B(\partdeg(l) = 0)\bigr)^k\\[-4pt]
\label{eq:minimum_nonbotnet_degree}
  &= 1 - \bigl((1 - p)^{n-k}\bigr)^k
  \geq 1 - \e^{-(1+\smallO(1)) n p k}
  \to 1 \,.
\end{align}
Hence, using \eqref{eq:nonbotnet_diameter_bound} and \eqref{eq:minimum_nonbotnet_degree}, we obtain
\begin{equation}
\label{eq:diameter_bound}
\P_B\Bigl(\max_{i, j \in V} \distG(i, j) \leq \bigO(1) / r\Bigr)
  = \P_B\bigl(\diam(G) \leq \bigO(1) / r\bigr)
  \to 1 \,.
\end{equation}

Finally, it follows from the dominated convergence theorem and \eqref{eq:distance_bound} that $\E_B[\1{\{\text{diam}(G) \leq \bigO(1) / r\}} \mspace{2mu} \distG^{\textup{avg}}(G)] = \smallO(1) / r$. Combining this with \eqref{eq:diameter_bound} and Markov's inequality we obtain, for any $a > 0$,
\begin{align}
\P_B\left(\distG^{\textup{avg}}(G) \geq \frac{a}{r}\right)
  &= \P_B\left(\1{\{\text{diam}(G) \leq \bigO(1) / r\}} \mspace{2mu} \distG^{\textup{avg}}(G) \geq \frac{a}{r}\right) - \smallO(1)\\
  &\leq \frac{r}{a} \, \E_B\bigl[\1{\{\text{diam}(G) \leq \bigO(1) / r\}} \mspace{2mu} \distG^{\textup{avg}}(G)\bigr] - \smallO(1)
  \to 0 \,.
\end{align}
In particular, choosing $a = (1 - \epsilon) \frac{d}{2(d+1)}$ gives $\P_B\bigl(\distG^{\textup{avg}}(G) < (1 - \epsilon) \frac{d}{2(d+1)}\frac{1}{r}\bigr) \to 1$. This shows that the average distance test is asymptotically powerful.
\end{proof}

%%%%%%%%%%%%%%%%%%%%%%%%%%%%%%%%%%%%%%%%%%%%%%%%%%%%%%%%%%%%%%%%%%%%%%%%%%%%%%%%
%%%%%%%%%%%%%%%%%%%%%%%%%%%%%%%%%%%%%%%%%%%%%%%%%%%%%%%%%%%%%%%%%%%%%%%%%%%%%%%%
\subsection{Proof of Theorem~\ref{thm:isolated_star_estimator_powerful}: Isolated star estimator performance}
\label{sec:proof_isolated_star_estimator_powerful}
\begin{proof}[\unskip\nopunct]
We need to show that $\estrisk(\est) \to 0$, for the estimator $\est$ from Definition~\ref{def:isolated_star_estimator}. First, we decompose the risk $\estrisk(\est)$ as
\begin{align}\label{eq:estimation_risk_decomposition}
\estrisk(\est)
  &= \E_B\left[\frac{|\est \symdiff B|}{2|B|}\right]
  = \frac{\E_B\bigl[| (V \setminus \est) \cap B|\bigr] + \E_B\bigl[|\est \cap (V \setminus B)|\bigr]}{2|B|}\\
  &= \frac{1}{2|B|} \sum_{j \in B} \P_B\bigl(j \notin \est\bigr) + \frac{1}{2|B|} \sum_{j \in V \setminus B} \P_B\bigl(j \in \est\bigr)\\
  &= \frac{1}{2|B|} \sum_{j \in B} \P_B\bigl(|S(j)| \leq \kissingnumber{d} + \xi_n\bigr) + \frac{1}{2|B|} \sum_{j \in V \setminus B} \P_B\bigl(|S(j)| > \kissingnumber{d} + \xi_n\bigr) \,.\;\;
\end{align}

We start by showing that the second term in \eqref{eq:estimation_risk_decomposition} vanishes. Note that, for any non-botnet vertex $i \in V \setminus B$, the size of its isolated star $|S(i)|$ is bounded by the kissing number $\kissingnumber{d}$ plus the amount of botnet vertices connected to it. Therefore,% for $i \in V \setminus B$, we obtain
\begin{multline}
\frac{1}{2|B|} \sum_{j \in V\setminus B} \P_B\bigl(|S(j)| > \kissingnumber{d} + \xi_n\bigr)
  = \frac{n - k}{2k} \, \P_B\bigl(|S(i)| > \kissingnumber{d} + \xi_n\bigr)\\[-3pt]
\begin{aligned}[b]
  &\leq \frac{n - k}{2k} \, \P_B\bigl(\text{$i$ is connected to at least $\xi_n$ botnet vertices}\bigr)\\
  &= \frac{n - k}{2k} \, \P\bigl(\text{Bin}(k, p) > \xi_n\bigr)
  \leq \frac{n - k}{2k} \left(\frac{k p \e}{\xi_n}\right)^{\xi_n}
  \to 0 \,,
\end{aligned}
\end{multline}
where the convergence to $0$ follows from the definition of $\xi_n$ in \eqref{eq:isolated_star_estimator_threshold_increase}. In fact, the definition of $\xi_n$ was chosen precisely to ensure this convergence.

To complete the proof, we analyze the first term on the right-hand side of \eqref{eq:estimation_risk_decomposition}. Let $i \in B$ be an arbitrary botnet vertex, then
\begin{multline}
\frac{1}{|B|} \sum_{j \in B} \P_B\bigl(|S(j)| \leq \kissingnumber{d} + \xi_n\bigr)
  = \P_B\bigl(|S(i)| \leq \kissingnumber{d} + \xi_n\bigr)\\[-3pt]
\begin{aligned}[b]
  &= 1 - \P_B\bigl(|S(i)| > \kissingnumber{d} + \xi_n\bigr)\\[3pt]
  &= 1 - \P_B\bigl(|S(i)| > \kissingnumber{d} + \xi_n \,\big|\, \text{deg}(i) > \kissingnumber{d} + \xi_n \bigr) \, \P_B\bigl(\text{deg}(i) > \kissingnumber{d} + \xi_n \bigr) \,.\notag
\end{aligned}
\end{multline}
Now, using the same argument as in \eqref{eq:isolated_star_probability_lower_bound}, we obtain
\begin{multline}
\P_B\bigl(|S(i)| > \kissingnumber{d} + \xi_n \,\big|\, \text{deg}(i) > \kissingnumber{d} + \xi_n \bigr)\\
  \geq \!\prod_{j=1}^{\kissingnumber{d}+\xi_n+1}\! \min\bigl\{\bigl(1 - (\kissingnumber{d} + \xi_n) p\bigr),(1-p)^{\kissingnumber{d} + \xi_n}\bigr\} \,,
\end{multline}%
which converges to $1$ provided that $\xi_n^2 p \to 0$. Combining the above, we obtain
\begin{align}
\label{eq:estrisk_bound}
\estrisk(\est)
  &= \frac{1}{2|B|} \sum_{j \in B} \P_B\bigl(|S(j)| \leq \kissingnumber{d} + \xi_n\bigr) + \frac{1}{2|B|} \sum_{j \in V \setminus B} \P_B\bigl(|S(j)| > \kissingnumber{d} + \xi_n\bigr)\\
  &= \frac12\bigl(1 - \bigl(1 - (\kissingnumber{d} + \xi_n) p\bigr)^{\kissingnumber{d} + \xi_n + 1} )\P_B\bigl(\text{deg}(i) > \kissingnumber{d} + \xi_n\bigr)\bigr) + \smallO(1) \,,
\end{align}
where $i \in B$ is an arbitrary botnet vertex. Therefore, the isolated star estimator has exact recovery when $\xi_n^2 p \to 0$  and $\P_B\bigl(\text{deg}(i) > \kissingnumber{d} + \xi_n\bigr) \to 1$, and partial recovery when $\xi_n^2 p \to 0$ and $\P_B\bigl(\text{deg}(i) > \kissingnumber{d} + \xi_n \bigr) = \bigOmega(1)$. To show this, we consider the three different cases from the theorem statement.

\paragraph{\normalfont\emph{Case \ref{thm:isolated_star_estimator_powerful_part1}:}}
From our assumption it follows that $k p \leq n^{-\alpha}$ for some $\alpha \in (0,1)$. Recall that $\lambertW_0(x)$ denotes the Lambert-W function, which can be approximated by $\lambertW_0(x) \asymp \log(x)$ when $x \to \infty$ \cite{Corless1996}. We obtain $\xi_n \asymp 2 \log(n / k) / \log(n^{\alpha}) = \bigO(1)$. Hence, it follows that $\xi_n^2 p \to 0$. Moreover,
\begin{align}
\P_B\bigl(\text{deg}(i) > \kissingnumber{d} + \xi_n \bigr)
  &= \P\bigl(\text{Bin}(n-1, p) > \bigO(1)\bigr)\\
  &= \begin{cases}
    1 - \smallO(1) & \text{if } n p \to \infty \,,\\
    \bigOmega(1)   & \text{otherwise} \,.
  \end{cases}
\end{align}
Therefore, the isolated star estimator achieves exact recovery when $n p \to \infty$, and partial recovery otherwise.

\paragraph{\normalfont\emph{Case \ref{thm:isolated_star_estimator_powerful_part2}:}}
From our assumption it follows that $n^{-\smallO(1)} \leq k p \leq \smallO(\log(n / k))$. Using $\lambertW_0(x) \to \infty$ when $x \to \infty$, we obtain
\begin{equation}
\xi_n \leq \frac{2 \log(n / k)}{\lambertW_0\bigl(\log(n / k) / \smallO(\log(n / k))\bigr)}
  = \smallO(\log(n / k)) \,.
\end{equation}
Hence, it follows that $\xi_n^2 p \leq \smallO(\log(n / k)^2) \log(n / k)^{-2} \to 0$. Moreover, from the assumptions for this case it follows that $n p \gg \log(n / k) \to \infty$, and therefore
\begin{align}
\P_B\bigl(\text{deg}(i) > \kissingnumber{d} + \xi_n \bigr)
  &= \P\bigl(\text{Bin}(n-1, p) > \kissingnumber{d} + \smallO(\log(n / k))\bigr)\\
  &\geq \P\bigl(\text{Bin}(n-1, p) > \log(n / k)\bigr)
  \to 1 \,.
\end{align}
Hence, the isolated star estimator has exact recovery.

\paragraph{\normalfont\emph{Case \ref{thm:isolated_star_estimator_powerful_part3}:}}
From our assumption it follows that $k p \geq \bigOmega(\log(n / k))$. When $k p \gg \log(n / k)$ we use that $\lambertW_0(x) \asymp x$ when $x \to 0$ \cite{Corless1996}, and obtain $\xi_n = \bigTheta(k p)$. Otherwise, when $k p = \bigTheta(\log(n / k))$, it also holds that $\xi_n = \bigTheta(\log(n / k)) = \bigTheta(k p)$. In both cases, it follows that $\xi_n^2 p = \bigTheta(k^2 p^3) \to 0$. Furthermore, note that $n p \gg k p \to \infty$ and therefore $\P_B\bigl(\text{deg}(i) > \kissingnumber{d} + \xi_n \bigr) = \P\bigl(\text{Bin}(n-1, p) > \bigO(k p)\bigr) \to 1$, so the isolated star estimator achieves exact recovery.
\end{proof}

%%%%%%%%%%%%%%%%%%%%%%%%%%%%%%%%%%%%%%%%%%%%%%%%%%%%%%%%%%%%%%%%%%%%%%%%%%%%%%%%
%%%%%%%%%%%%%%%%%%%%%%%%%%%%%%%%%%%%%%%%%%%%%%%%%%%%%%%%%%%%%%%%%%%%%%%%%%%%%%%%
\subsection{Proof of Theorem~\ref{thm:no_test_powerful}: No test is powerful}
\label{sec:proof_no_test_powerful}
\begin{proof}[\unskip\nopunct]
We start by considering a simpler version of the problem where the set of potential botnet vertices $B \subseteq V$ is known. Now, we no longer have a composite alternative hypothesis, and this problem corresponds to a hypothesis test between two simple hypotheses. That is, given a set $B \subseteq V$, we consider the risk
\begin{equation}
\label{eq:worst_case_risk_simple_vs_simple}
\detrisk^{*}(\test) = \P_0(\test(G) \neq 0) + \P_B(\test(G) \neq 1) \,.
\end{equation}
Note that, for every test $\test$, the risk $\detrisk^{*}(\test)$ is a lower bound for the worst-case risk $\detrisk(\test)$ in \eqref{eq:worst_case_risk}. Using a result by Tsybakov \cite[Proposition 2.1]{Tsybakov2009}, for every test $\test$ it follows that
\begin{equation}
\label{eq:worst_case_risk_simple_vs_simple_lower_bound}
\detrisk(\test)
  \geq \detrisk^{*}(\test)
  \geq \sup_{\tau > 0} \, \left\{ \frac{\tau}{\tau + 1} \, \P_0(L(G) \geq \tau)\right\} \,,
\end{equation}
where $L(g) = \P_B(G = g) / \P_0(G = g)$ is the likelihood ratio. Therefore, to show that no test is asymptotically powerful it suffices to show that $\P_0(L(G) \geq \tau)$ remains bounded away from zero, for some $\tau$ independent of the graph size $n$. To this end, define the event
\begin{equation}
\label{eq:isolated_event}
\eventib % Event isolated botnet
  \coloneqq \bigl\{\text{all vertices in $B$ are isolated in the graph $G$}\bigr\} \,.
\end{equation}
For every graph $g$ such that $\P_0(G = g\,|\,\eventib) > 0$ (i.e., a graph that could be a sample from the null hypothesis with all vertices in $B$ being isolated), it follows that
\begin{align}
\P_0(G = g)
  \leq \P_0(G_{\scriptscriptstyle V \setminus B} = g_{\scriptscriptstyle V \setminus B})
  = \P_B(G_{\scriptscriptstyle V \setminus B} = g_{\scriptscriptstyle V \setminus B})
  = \frac{\P_B(G = g)}{(1 - p)^{(n - k)k + k(k-1)/2}} \,,
\end{align}
where we have used $\{G_{\scriptscriptstyle V \setminus B} = g_{\scriptscriptstyle V \setminus B}\}$ to indicate the event where the subgraphs induced by the non-botnet vertices $V \setminus B$ are equal. Hence, for all $g$ in which the vertices of $B$ are isolated, we obtain
\begin{equation}
\label{eq:likelihood_ratio_bound2}
L(g)
  = \frac{\P_B(G = g)}{\P_0(G = g)}
  \geq (1 - p)^{(n - k)k + k(k-1)/2}
  = \e^{-(1 + \smallO(1)) n p k} \,,
\end{equation}
which remains strictly positive as $n \to \infty$ by the assumption that $n p k = \bigO(1)$.\pagebreak[3]
Therefore, we can choose $\tau > 0$ small enough such that $\P_0(L(G) \geq \tau \,|\, \eventib) = 1$ for all $n$ large enough. Finally, using the same reasoning as in \eqref{eq:isolated_star_probability_lower_bound}, observe that
\begin{align}
\P_0(L(G) \geq \tau)
  &\geq \P_0(L(G) \geq \tau \,|\, \eventib) \, \P_0(\eventib)\\
  &= \P_0(\eventib)\\
  &\geq \left(\prod_{i=0}^{k-1} (1 - i p)\right) (1 - k p)^{n - k}\\
  &= \e^{-(1 + \smallO(1)) n p k} \,.
\end{align}
which remains strictly positive as $n \to \infty$ by the assumption that $n p k = \bigO(1)$. Plugging this into \eqref{eq:worst_case_risk_simple_vs_simple_lower_bound} shows that, for every test $\test$, the risk $\detrisk(\test) \geq \detrisk^{*}(\test)$ remains bounded away from zero, and therefore that no test can be asymptotically powerful.
\end{proof}

%%%%%%%%%%%%%%%%%%%%%%%%%%%%%%%%%%%%%%%%%%%%%%%%%%%%%%%%%%%%%%%%%%%%%%%%%%%%%%%%
%%%%%%%%%%%%%%%%%%%%%%%%%%%%%%%%%%%%%%%%%%%%%%%%%%%%%%%%%%%%%%%%%%%%%%%%%%%%%%%%
\subsection{Proof of Lemma~\ref{lem:dimension_estimator_consistent}: Dimension estimator is consistent}
\label{sec:proof_dimension_estimator_consistent}
\begin{proof}[\unskip\nopunct]
We start by showing that $\widehat{C}_d \conv[\P_0] C_d$, and from this it follows that $\widehat{d} \conv[\P_0] d$ by the continuous mapping theorem and because \eqref{eq:clustering_coefficient_analytic} is continuous. Using \eqref{eq:clustering_coefficient_estimate} we obtain
\begin{align}
\label{eq:clustering_coefficient_estimate2}
\widehat{C}_d(G)
  &= \frac{n^{-3} \sum_{1 \leq i, j, k \leq n} \1{\{i \leftrightarrow j, i \leftrightarrow k, j \leftrightarrow k\}} / p^2}{n^{-3} \sum_{1 \leq i, j, k \leq n} \1{\{i \leftrightarrow j, i \leftrightarrow k\}} / p^2} \,.
\end{align}
Here we will show that the numerator converges in probability to $C_d$, and the denominator converges in probability to $1$. Since the computations regarding the denominator are largely similar to those of the numerator these will be omitted for brevity, and we will focus on the numerator.

Let $X_{ijk} = \1{\{i \leftrightarrow j, i \leftrightarrow k, j \leftrightarrow k\}} / p^2$ and $\bar{X} = n^{-3} \smash{\sum_{1 \leq i, j, k \leq n}} X_{ijk}$, then $\bar{X}$ is precisely the numerator in \eqref{eq:clustering_coefficient_estimate2}. Consider the first moment of $\bar{X}$, this is given by
\begin{align}
\hspace{10pt}&\hspace{-10pt}
\E_0[\bar{X}]
  = n^{-3} \hspace{1.2em}\smashoperator{\sum_{1 \leq i, j, k \leq n}}\hspace{1.2em} \E_0[X_{ijk}]\\
  &= n^{-3} \hspace{1.2em}\smashoperator{\sum_{1 \leq i, j, k \leq n}}\hspace{1.2em} \frac{\P_0(i \leftrightarrow j)\P_0(i \leftrightarrow k)\P_0(j \leftrightarrow k \,|\, i \leftrightarrow j, i \leftrightarrow k)}{p^2}
  = (1 + \smallO(1)) \, C_d \,. \notag
\end{align}
Moreover, the second moment of $\bar{X}$ can be computed by splitting between the number of common vertices in the two triangles involved. This gives
\begin{align}
\E_0[\bar{X}^2]
  &= n^{-6} \hspace{1.2em}\smashoperator{\sum_{1 \leq i, j, k, i', j', k' \leq n}}\hspace{1.2em}
    \E_0[X_{ijk} \, X_{i'j'k'}]\\
  &= n^{-6} \hspace{1.2em}\smashoperator{\sum_{\substack{1 \leq i, j, k, i', j', k' \leq n\\\text{distinct}}}}\hspace{1.2em}
    \E_0[X_{ijk} \, X_{i'j'k'}]
   \;+\; 3 \, n^{-6} \hspace{1.2em}\smashoperator{\sum_{\substack{1 \leq i, j, k, j', k' \leq n\\\text{distinct}}}}\hspace{1.2em}
     \E_0[X_{ijk} \, X_{ij'k'}]\\
  &\qquad\qquad \;+\; 3 \, n^{-6} \hspace{1.2em}\smashoperator{\sum_{\substack{1 \leq i, j, k, k' \leq n\\\text{distinct}}}}\hspace{1.2em}
    \E_0[X_{ijk} \, X_{ijk'}]
    \;+\; n^{-6} \hspace{1.2em}\smashoperator{\sum_{\substack{1 \leq i, j, k \leq n\\\text{distinct}}}}\hspace{1.2em}
    \E_0[X_{ijk}^2]\\
%  &= n^{-6} \sum_{\substack{1 \leq i, j, k, i', j', k' \leq n\\\text{distinct}}} C_d^2 + 3 \, n^{-6} \sum_{\substack{1 \leq i, j, k, j', k' \leq n\\\text{distinct}}} C_d^2\\
%  &\qquad\qquad + 3 \, n^{-6} \sum_{\substack{1 \leq i, j, k, k' \leq n\\\text{distinct}}} C_d^2 / p + n^{-6} \sum_{\substack{1 \leq i, j, k \leq n\\\text{distinct}}} C_d / p^2\\
  &= (1 + \smallO(1)) \left[C_d^2 + 3 \, \frac{C_d^2}{n} + 3 \, \frac{C_d^2}{n^2 p} + \frac{C_d}{n^3 p^2}\right]
  = (1 + \smallO(1)) \, C_d^2 \,,
\end{align}
where the final step follows from the assumption that $p \geq \bigOmega(1 / n)$. Hence, $\text{Var}_0(\bar{X}) = \E_0[\bar{X}^2] - \E_0[\bar{X}]^2 = \smallO(1)$, and therefore it follows by Chebyshev's inequality that $\bar{X} \conv[\P_0] C_d$. This shows that the numerator of \eqref{eq:clustering_coefficient_estimate2} converges in probability to $C_d$ and the denominator of \eqref{eq:clustering_coefficient_estimate2} converges in probability to $1$, so we have $\widehat{C}_d \conv[\P_0] C_d$. Finally, it follows from the continuous mapping theorem that $\widehat{d} \conv[\P_0] d$, and we conclude that our estimator for the dimension is consistent under the null hypothesis.

Under the alternative hypothesis, the proof is largely similar. Because the botnet size $k = \smallO(n)$ is small, it can be seen that the first and second moment of $\bar{X}$ converge to the same values, and therefore $\bar{X} \conv[\P_B] C_d$. Finally, we can again apply the continuous mapping theorem to show that our estimator for the dimension is consistent under the alternative hypothesis.
\end{proof}

%%%%%%%%%%%%%%%%%%%%%%%%%%%%%%%%%%%%%%%%%%%%%%%%%%%%%%%%%%%%%%%%%%%%%%%%%%%%%%%%
%%%%%%%%%%%%%%%%%%%%%%%%%%%%%%%%%%%%%%%%%%%%%%%%%%%%%%%%%%%%%%%%%%%%%%%%%%%%%%%%
\subsection{Proof of Lemma~\ref{lem:probability_estimator_consistent}: Connection probability estimator is consistent}
\label{sec:proof_connection_probability_estimator_consistent}
\begin{proof}[\unskip\nopunct]
We start by showing that $\widehat{p} / p \conv[\P_0] 1$. Using the estimator $\widehat{p}$ from \eqref{eq:edge_probability_estimate} it follows directly that $\E_0\bigl[\widehat{p} / p\bigr] = 1$. Therefore, we are left to compute
\begin{align}
\E_0[(\widehat{p} / p)^2]
  &= \binom{n}{2}^{\!\!-2} \!\! \sum_{\substack{\cramped{1\leq i<j\leq n}\\\cramped{1\leq i' < j' \leq n}}}\E_0\left[\frac{\1{\{i \leftrightarrow j\}}}{p}\frac{\1{\{i' \leftrightarrow j'\}}}{p}\right] \\
  &= \binom{n}{2}^{\!\!-2} \biggl(\biggl[\binom{n}{2}^{\!2} - \binom{n}{2}\biggr] + \binom{n}{2}\frac{1}{p}\biggr)\\
  &= \biggl(1 - \binom{n}{2}^{\!\!-1} + \binom{n}{2}^{\!\!-1} \, \frac{1}{p}\biggr)
  = 1 + \smallO(1) \,,
\end{align}
where we obtained the second equality by splitting between the case where $i \neq i'$ and $j \neq j'$, and the case where $i = i'$ and $j = j'$. Moreover, the final step followed from the assumption $p \geq \bigOmega(1/n)$. 
Therefore, it follows that $\text{Var}_0(\widehat{p} / p) = \E_0\bigl[(\widehat{p} / p)^2\bigr] - \E_0\bigl[\widehat{p} / p\bigr]^2 = \smallO(1)$, and hence $\widehat{p} / p \conv[\P_0] 1$ by Chebyshev's inequality. Moreover, for any distinct triplet $i, j, k \in V$,
\begin{equation}
p   = \P_0(i \leftrightarrow j) = \P_B(i \leftrightarrow j) \,, \qquad
p^2 = \P_0(i \leftrightarrow j, i \leftrightarrow k) = \P_B(i \leftrightarrow j, i \leftrightarrow k) \,.
\end{equation}
Hence, performing the above computations under the measure $\P_B$ shows that $\widehat{p}/p \conv[\P_B] 1$ as well. 
\end{proof}

\paragraph{Acknowledgements.} GB acknowledges the support of the STAR cluster and Eurandom for visiting KB, RC, and RvdH at TU/e. The work of RvdH is supported by the NWO Gravitation Networks grant 024.002.003.

%%%%%%%%%%%%%%%%%%%%%%%%%%%%%%%%%%%%%%%%%%%%%%%%%%%%%%%%%%%%%%%%%%%%%%%%%%%%%%%%
%%%%%%%%%%%%%%%%%%%%%%%%%%%%%%%%%%%%%%%%%%%%%%%%%%%%%%%%%%%%%%%%%%%%%%%%%%%%%%%%
\printbibliography[title={References}]
\end{document}